\documentclass[11pt,
]{article}
\usepackage{amssymb,amstext,amsmath,amsthm,latexsym}
\usepackage{mathrsfs,mathbbol}
\newcommand{\nek}{\newcommand}
\nek{\renek}{\renewcommand}

\DeclareMathAlphabet{\cur}{U}{eur}{m}{n}
\DeclareMathAlphabet{\cub}{U}{eur}{b}{n}
\nek{\skr}{\mathscr}

3

\nek{\bfsf} {\sffamily\bfseries\upshape}

\nek{\bfit}{\bfseries\itshape}
\nek{\bfsl}{\bfseries\slshape}
\nek{\sfbs}{\mdseries\sffamily\itshape}        %\bfseries
\DeclareMathAlphabet{\bma}{OML}{cmm}{b}{it}

\nek{\vyk} [1] {}

\nek{\imar} [1] {}

\vyk{
\nek{\imar}[1]{\marginpar[%\vspace{-1ex}%
\flushright\footnotesize\sl
$\mtho\longrightarrow$\\ \vspace{-1ex}{#1}$\mtho\dashv$\vspace*{1ex}]
{%\vspace{-1ex}%
\flushleft\footnotesize\sl
$\mtho\longleftarrow$\\ \vspace{-1ex}{#1}$\mtho\dashv$\vspace*{1ex}}}

\nek{\jmar}[1]{\marginpar[\vspace{5ex}%
\flushright\footnotesize\sl
$\mtho\longrightarrow$\\ \vspace{0ex}{#1}$\mtho\dashv$\vspace*{0ex}]
{\vspace{5ex}%
\flushleft\footnotesize\sl
$\mtho\longleftarrow$\\ \vspace{0ex}{#1}$\mtho\dashv$\vspace*{0ex}}}
}

\renek{\thesubsection}{\arabic{subsection}}

\nek{\punk}[1]{\subsection{{\protect\boldmath#1}}%
\setcounter{clt}0}

\renek{\thefigure}{\arabic{figure}}

\renek{\subsectionmark}[1]{\markright{%
{\small\upshape\thesubsection\hspace{1em}#1}}
{}%
}%

\newcounter{enuf}
\nek{\enufi}{\addtocounter{enuf}{1}}
\nek{\fenu}{

\enufi\itsep
}
\newcounter{enuc}
\nek{\enuci}{\addtocounter{enuc}{1}}
\nek{\cenu}{

\enuci\itsep
}

\newcounter{enuF}
\nek{\enuFi}{\addtocounter{enuF}{1}}
\nek{\Fenu}{

\nek{\ifla}{\addtocounter{enuF}1\itla}
\enuFi\itsep
}

\nek{\itsep}{\itemsep=0.4ex plus 0.15ex minus 0.15ex}
\nek{\tenu}[1]{

\itsep
}
\nek{\tenui}[1]{
\def\theenumii{#1}
\itsep
}

\theoremstyle{plain}

\newtheorem{theore}      {Theorem}               
\newtheorem{corollar}  [theore]{Corollary}
\newtheorem{propo}     [theore]{Proposition}
\newtheorem{lemm}      [theore]{Lemma}
\newtheorem{cla}       [theore]{Claim} 
\newtheorem{clt}       {Claim} [theore]
\newtheorem{cort}  [clt]{Corollary}  %[theore]
\newtheorem{prot}  [clt]{Proposition}  %[theore]
\newtheorem{lemt}  [clt]{Lemma}

\newtheorem{cly}       {{\ubf Claim}} 
\newtheorem{cory}    [cly]  {{\ubf Corollary}}  %[theore]

\theoremstyle{definition}
\newtheorem{defn}      [theore]{Definition}
\newtheorem{prim}      [theore]{Example}
\newtheorem{primt}     {Example}[theore]
        \renek{\theprimt}{\thetheore:\roman{primt}}
        \nek{\sbrospri}{\addtocounter{theore}1\setcounter{primt}0}
\newtheorem{rem}       [theore]{Remark}
\newtheorem{que}       [theore] {Question}  
\newtheorem{remt}  [clt]{Remark}  %[theore]
\newtheorem{vop}       [theore]{Question}
\newtheorem{upt}       [theore]{Exercise}
\newtheorem{uprt}       [clt]{Exercise}
\newtheorem{deft}  [clt]{Definition}

\newtheorem*{qrF}{{\ubf Proof}}               %\def\theprF{}
\newtheorem*{prF}{{\bf Proof}}               %\def\theprF{}
\newtheorem*{prC}{{\bf Proof of the claim}}  %\def\theprC{}
\newtheorem*{con}{{\bfit Construction\/}{\bf.}}    %\def\thecon{}
\newtheorem*{aq} {{\bf Acknowledgements\/}}          %\def\theaq{}

\nek{\thsp}{\hspace{0.1ex plus \mathsurround}}
\nek{\bpro}{\begin{propo}}
\nek{\epro}{\end{propo}}
\nek{\bcl} {\begin{cla}}
\nek{\ecl} {\end{cla}}
\nek{\bct} {\begin{clt}}
\nek{\ect} {\end{clt}}
\nek{\bcy} {\begin{cly}}
\nek{\ecy} {\end{cly}}
\nek{\bco} {\begin{con}}
\nek{\eco} {\qeDD{Construction}\end{con}}
\nek{\bcor}{\begin{corollar}}
\nek{\ecor}{\end{corollar}}
\nek{\bcoy}{\begin{cory}}
\nek{\ecoy}{\end{cory}}
\nek{\baq} {\begin{aq}}
\nek{\eaq} {\end{aq}}
\nek{\bex} {\begin{prim}}
\nek{\eex} {\qeD\end{prim}}
\nek{\eeX} {\end{prim}}
\nek{\bpri} {\begin{primt}}
\nek{\epri} {\qeD\end{primt}}
\nek{\bdf} {\begin{defn}} %\rm\thsp}
\nek{\eDf} {\end{defn}}
\nek{\edf} {\qeD\end{defn}}
\nek{\edF} {\end{defn}}
\nek{\ble} {\begin{lemm}}
\nek{\ele} {\end{lemm}}
\nek{\bte} {\begin{theore}}
\nek{\ete} {\end{theore}}
\nek{\bre} {\begin{rem}} 
\nek{\ere} {\qeD\end{rem}} 
\nek{\bqe} {\begin{que}} 
\nek{\eqe} {\qeD\end{que}}
\nek{\bdt}{\begin{deft}\thsp\rm}
\nek{\edt}{\qeD\end{deft}}
\nek{\blt}{\begin{lemt}}
\nek{\elt}{\end{lemt}}
\nek{\bcot}{\begin{cort}}
\nek{\ecot}{\end{cort}}
\nek{\bprt}{\begin{prot}}
\nek{\eprt}{\end{prot}}
\nek{\bup} {\begin{upt}}
\nek{\eup} {\qeD\end{upt}}
\nek{\eUp} {\end{upt}}
\nek{\bqu} {\begin{vop}}
\nek{\equ} {\qeD\end{vop}}
\nek{\bret}{\begin{remt}}
\nek{\eret}{\qeD\end{remt}}
\nek{\bupt} {\begin{uprt}}
\nek{\eupt} {\qeD\end{uprt}}
\nek{\bqf} {\begin{qrF}} 
\nek{\eqf} {\,\hfill{\mtho$\vdash$}\end{qrF}} 

\nek{\bpf} {\begin{prF}} 
\nek{\epf} {\qed\end{prF}} 
\nek{\ePf} {\end{prF}} 
\nek{\bpc} {\begin{prC}} 
\nek{\epc} {\qeDD{Claim}\end{prC}} 
\nek{\qeD} {\hfill$\mtho\Box$}
\nek{\qfD} {\hfill$\mtho\vdash$}
\nek{\qeDD} [1] 
{\hfill\hbox{$\mtho\Box$~({\small\sl #1\/}\hspace{0.1ex})}}
\nek{\qedd} [1] 
{\hfill\hbox{$\mtho\vdash$~({\sl #1\/}\hspace{0.1ex})}}
\nek{\epF} [1] {\qeDD{#1}\end{prF}} 
\nek{\eqF} [1] {\qedd{#1}\end{qrF}} 

\nek{\bde}{\begin{description}}
\nek{\ede}{\end{description}}
\nek{\ben}{\begin{enumerate}}
\nek{\een}{\end{enumerate}}
\nek{\bit}{\begin{itemize}}
\nek{\eit}{\end{itemize}}
\nek{\bay}{\begin{array}}
\nek{\eay}{\end{array}}
\nek{\bmp}{\begin{minipage}}
\nek{\emp}{\end{minipage}}
\nek{\bus}{\begin{equation}}   
\nek{\eus}{\end{equation}}

\nek{\btb}{\begin{tabular*}}
\nek{\etb}{\end{tabular*}}
\nek{\fF} {{\bf F}}
\nek{\fG} {{\bf G}}
\nek{\Gd} {\fG_\fda}
\nek{\Gds}{\fG_{\fda\fsg}}
\nek{\Fs} {\fF_\fsg}
\nek{\fS} {{\bf S}}
\nek{\fT} {{\bf T}}
\nek{\fH} {{\bf H}}
\nek{\fa} {{\bf a}}
\nek{\fb} {{\bf b}}
\nek{\ZFC}{\text{\ubf ZFC}}
\nek{\ZC} {\text{\ubf ZC}}
\nek{\zhc}{\text{\ubf ZFC}^-}
\nek{\pli}{{\bf I}}
\nek{\plj}{{\phantom{I}\bf I}}
\nek{\pld}{{\bf II}}
\nek{\gai} [2] {\fH^\cZ_{#1#2}}
\nek{\gad} [2] {\fG^{#1}_{#2}} %{\fG^\cZ_{#1#2}}
\nek{\iai} [1] {\fH_{#1}}
\nek{\iad} [1] {\fG_{#1}}
\nek{\hi} {\fH}
\nek{\hd} {\fG}
\nek{\etc} {{\sl etc}}
\nek{\iesp}{\hspace{0.3ex}}
\nek{\pw} {\hbox{a.\iesp e.}}
\nek{\ie} {\hbox{\sl i.\iesp e.}}
\nek{\eg} {\hbox{\sl e.\iesp g.}}
\nek{\ea} {\hbox{\sl et.\hspace{0.3ex}al.}}
\nek{\noo}
{\text{\sl w.\iesp l.\iesp o.\iesp g.}}
\nek{\pv} {\text{a.\iesp a.}}
\nek{\vrt} {\text{w.\iesp r.\iesp t.}}
\nek{\vva} {{\sl vice versa}}
\nek{\lsc} {\text{\sc lsc}}
\nek{\er} {{ER}}
\nek{\bm} {{BM}}
\nek{\dd}[1]{$\mtho\hspace{0.2ex}{#1}$-\hspace{0.0ex}}
\nek{\dw}{\dd\om}
\nek{\lis}[1] {\mathop{\tt lim\hspace{0.2ex}sup}_{#1}}
\nek{\len}[1] {\mathop{\tt lh}{#1}}
\nek{\Ord}  {{\tt{Ord}}}
\nek{\Exh}  {{\tt{Exh}}}
\nek{\Nul}  {{\tt{Null}}}
\nek{\Mod}  {\mathop{\tt{Mod}}}
\nek{\Aut}  {\mathop{\tt{Aut}}}
\nek{\card} {\mathop{\tt card}}
\nek{\lh}   {\mathop{\tt lh}}
\nek{\pr}   {\mathop{\tt pr}}
\nek{\sr}   {\mathop{\tt sr}}
\nek{\ran}  {\mathop{\tt ran}}
\nek{\dom}  {\mathop{\tt dom}}
\nek{\fld}  {\mathop{\tt field}}
\nek{\otp}  {\mathop{\tt otp}}
\nek{\Max}  {\mathop{\tt Max}}
\nek{\tsup} {\mathop{\tt sup}}
\nek{\tinf} {\mathop{\tt inf}}
\nek{\tmin} {\mathop{\tt min}}
\nek{\tmax} {\mathop{\tt max}}
\nek{\tlim} {\mathop{\tt lim}}
\nek{\tlis} {\mathop{\tt lim\hspace{0.3ex}sup}}
\nek{\tlii} {\mathop{\tt lim\hspace{0.3ex}inf}}
\nek{\Fin}  {{\tt Fin}} 
\nek{\bFin} {{\bf Fin}} 
\nek{\maxi}[1] {\Max^\xi_{#1}}
\nek{\HC} {{\rm HC}}
\nek{\ccc}{{\sc ccc}}
\nek{\al} {\alpha}
\nek{\ba} {\beta}
\nek{\ga} {\gamma}
\nek{\Da} {\Delta}
\nek{\da} {\delta}
\nek{\ka} {\kappa}
\nek{\la} {\lambda}
\nek{\La}{\Lambda}
\nek{\sg} {\sigma}
\nek{\Sg} {\Sigma}
\nek{\vpi}{\varphi}
\nek{\vpy}{\vpi_\iy}
\nek{\vt} {\vartheta}
\nek{\vT} {\Theta}
\nek{\ovt}{{\overline\vt}}
\nek{\ovi}{{\overline\vpi}}
\nek{\ops}{{\overline\psi}}
\nek{\ve} {\varepsilon}
\nek{\om} {\omega}
\nek{\Om} {\Omega}
\nek{\lom}{^{<\om}}
\nek{\za}{\zeta}
\nek{\tpi}{\tau_\vpi}
\nek{\omi} {\om_1}
\nek{\omm} [1] {\om^{\om^{#1}}}
\nek{\bse} {2\lom}
\nek{\nse} {\dN\lom}
\nek{\alo} {{\aleph_0}}
\nek{\sd}   {\mathbin{\Da}}
\nek{\bigd} {\mathbin{\hbox{\large$\mtho\Da$}}}
\nek{\sqe} {\fmu\hspace{0.05ex}}
\newcommand{\iSg}{{\mathchar"7106}}
\newcommand{\iPi}{{\mathchar"7105}}
\newcommand{\iDa}{{\mathchar"7101}}
\newcommand{\is}[2]{\iSg^{#1}_{#2}}
\newcommand{\ip}[2]{\iPi^{#1}_{#2}}
\newcommand{\id}[2]{\iDa^{#1}_{#2}}
\nek{\BBB}{\hspace{0.05ex}}
\nek{\dA}{{\BBB{\mathbb A}\BBB}}
\nek{\dC}{{\BBB{\mathbb C}\BBB}}
\nek{\dF}{{\BBB{\mathbb F}\BBB}}
\nek{\dN}{{\BBB{\mathbb N}\BBB}}
\nek{\dP}{{\BBB{\mathbb P}\BBB}}
\nek{\dQ}{{\BBB{\mathbb Q}\BBB}}
\nek{\dqp}{\dQ^+}
\nek{\dR}{{\BBB{\mathbb R}\BBB}}
\nek{\dS}{{\BBB{\mathbb S}\BBB}}
\nek{\dT}{{\BBB{\mathbb T}\BBB}}
\nek{\dV}{{\BBB{\mathbb V}\BBB}}
\nek{\dZ}{{\BBB{\mathbb Z}\BBB}}
\nek{\dX}{{\BBB{\mathbb X}\BBB}}
\nek{\dY}{{\BBB{\mathbb Y}\BBB}}
\nek{\dyn} {\dY^\dN}
\nek{\dvp} {\dV^+}
\nek{\dn}{2^\dN}
\nek{\dntn}{2^{\dN\ti\dN}}
\nek{\dnqn}{{(2^\dN)}{}^\dN}
\nek{\pnqn}{{\pn}{}^\dN}
\nek{\bn}{\dN^\dN}
\nek{\rn} {\dR^\dN}
\nek{\tm} [1] {2^{\mxi\xi}}
\nek{\nn}{{\dN\ti\dN}}
\nek{\ccs} {}
\nek{\cA}{{\ccs{\skr A}\ccs}}
\nek{\cD}{{\ccs{\skr D}\ccs}}
\nek{\cE}{{\ccs{\skr E}\ccs}}
\nek{\cF}{{\ccs{\skr F}\ccs}}
\nek{\cS}{{\ccs{\skr S}\ccs}}
\nek{\cP}{{\ccs{\skr P}\ccs}}
\nek{\cW}{{\ccs{\skr W}\ccs}}
\nek{\cX}{{\ccs{\skr X}\ccs}}
\nek{\cI} {{\skr I}} 
\nek{\cJ} {{\skr J}} 
\nek{\cO} {{\skr O}} 
\nek{\cZ} {{\skr Z}}
\nek{\zo} {\cZ_0}
\nek{\zw} {\cZ_{\hbox{\small\rm w}}}
\nek{\xn} {\cX^\dN}
\nek{\an} {\dA^\dN}
\nek{\cd} [1] {\cD_{#1}}
\nek{\pws}  [1] {\cP(#1)}
\nek{\cp}  [2] {\cP_{\hspace*{-0.4ex}\tt cnt}^{#1}(#2)}
\nek{\pwf} [1] {\cP_{\hspace*{-0.4ex}\tt fin}(#1)}
\nek{\pwc} [1] {\cP_{\hspace*{-0.4ex}\tt ctbl}(#1)}
\nek{\pnn}{\cP(\nn)}
\nek{\pn}{\cP(\dN)}
\nek{\ps}{\cP(\dS)}
\nek{\pz}{\cP(\dZ)}
\nek{\mm}{{\BBB{\mathfrak M}\BBB}}
\nek{\mn}{{\BBB{\mathfrak N}\BBB}}
\nek{\gE}{{\BBB{\mathfrak E}\BBB}}
\nek{\shi} {{\mathfrak s}}
\newcommand{\gc}  {{\BBB{\mathfrak c}\BBB}}

\nek{\ilo}[1] {{[0\hspace{0.1ex},\hspace{0.1ex}n_{#1})}}
\nek{\ii} [1] {{[n_{#1}\hspace{0.1ex},\hspace{0.1ex}\infty)}}
\nek{\iv} [2] {{(#1\hspace{0.1ex},\hspace{0.1ex}#2)}}
\nek{\ix} [2] {{[#1\hspace{0.1ex},\hspace{0.1ex}#2]}}
\nek{\ir} [2] {{[#1\hspace{0.1ex},\hspace{0.1ex}#2)}}
\nek{\iry}[1] {{[#1\hspace{0.1ex},\hspace{0.1ex}\iy)}}
\nek{\iq} [2] {\ir{\nu_{#1}}{\nu_{#2}}}
\nek{\iqy}[1] {\ir{\nu_{#1}}\iy}
\nek{\iqo}[1] {\ir0{\nu_{#1}}}
\nek{\iqn}[1] {\iq{#1}{#1+1}}
\nek{\opl} {\oplus}
\nek{\ap}  {\cdot}
\nek{\cj}  {\mathbin{\hspace{0.1ex}\&\hspace{0.1ex}}}
\nek{\dm}  {$$}
\nek{\sus} {\mathopen{\exists\hspace{0.35ex}}}
\nek{\kaz} {\mathopen{\forall\hspace{0.35ex}}}
\nek{\imp} {\Longrightarrow} 
\nek{\eqv} {\Longleftrightarrow} 
\nek{\ti}  {\times} 
\nek{\mo}  {\models} 
\nek{\sq}  {\subseteq}
\nek{\qs}  {\supseteq}
\nek{\su}  {\subset}
\nek{\sneq}{\subsetneqq}
\nek{\we}  {{\mathbin{\hspace*{0.2ex}^\wedge}}}
\nek{\obr} {^{-1}}
\nek{\dif} {\smallsetminus}
\nek{\res} {\mathbin{\restriction}}
\nek{\lef} {\preccurlyeq}
\nek{\gef} {\succcurlyeq}
\nek{\pu}  {\varnothing}
\nek{\iy}  {\infty}
\nek{\piy} {+\iy}
\nek{\nin} {\not\in}
\nek{\limp}{\,\imp\,}
\nek{\leqv}{\,\eqv\,}
\nek{\onto}{\stackrel{{\rm onto}}{\longrightarrow}}
\nek{\ang} [1] {\langle #1\rangle}
\nek{\stk} [2] {\ang{#1\hspace{0.3ex};\hspace{0.1ex}#2}}
\nek{\sis} [2] {\ans{#1}_{#2}}
\nek{\ans} [1] {\{\hspace{0.01ex}#1\hspace{0.01ex}\}}

\nek{\zz} {\linebreak[0]} 
\nek{\ens} [2] {\ans{{#1\hspace{0.5ex}{:}}\zz\hspace{0.5ex}#2}}

\nek{\suh} [1] {[\hspace{0.3pt}#1\hspace{0.3pt}]_{\sq}}
\nek{\itla} {\item\label}

\nek{\aeq} {\mathbin{\|}}
\nek{\laeq}{\,\aeq\,} 

\nek{\tz} {\mathbin{;}}

\nek{\seq} [2] {(#1)_{#2}}
\nek{\seqq}[3] {{#1}_{#2}[#3]}

\nek{\kb}[2]{#1^{(#2)}}

\nek{\ssty} {\textstyle}
\nek{\isum} [3] {{{\ssty\ugl\sum_{#1}\,#3}\mid{#2}\ugr}}

\nek{\ifi}  {{\cur{Fin}\hspace*{0.1ex}}}
\nek{\frt}  {{\cur{Fr}\hspace*{0.1ex}}}
\nek{\ibo}  {{\cur{Bou}\hspace*{0.1ex}}}
\nek{\fifi} {\ifi\ti\ifi}
\nek{\fio} {\ifi\ti0}
\nek{\ofi} {0\ti\ifi}

\nek{\xip} {{\xi+1}}
\nek{\etp} {{\eta+1}}

\nek{\vx} [1] {^{(#1)}}

\nek{\dop}[1] {{#1}^{\complement}}

\nek{\skl} {\hbox{\mtho\large$($}}
\nek{\skp} {\hbox{\mtho\large$)$}}

\nek{\ugl} {\hbox{\mtho\large$\langle$}}
\nek{\ugr} {\hbox{\mtho\large$\rangle$}}

\nek{\df} [1] {\dop {#1}}
\nek{\pl} [1] {{#1}^+}

\nek{\bbW} {\hbox{\mtho\boldmath $W$}}

\nek{\bo}{{\bf O}}

\nek{\Zo}  {{\cZ_0}}
\nek{\Eo}  {\rE_{\text{\sf0}}}
\nek{\Fo}  {\rF_0}
\nek{\ovli} [1] {\widehat{#1}}
\nek{\Ec}  {\mathbin{\ovli{{\rE}}}}
\nek{\Eco} {\mathbin{\ovli{\Eo}}}
\nek{\rzo}  {\rZ_{\text{\sf0}}}
\nek{\neo} {\mathbin{\not{{\hspace{-0.4ex}\rE}}_0}}
\nek{\nei} {\mathbin{\not{{\hspace{-0.4ex}\rE}}_1}}
\nek{\nfo} {\mathbin{\not{{\hspace{-0.4ex}\rF}}_0}}
\nek{\nrE} [1] {\mathbin{\not{{\hspace{-0.4ex}\rE}}_{#1}}}
\nek{\rC}  {\mathbin{\sf C}}
\nek{\rD}  {\mathbin{\sf D}}
\nek{\rR}  {\mathbin{\sf R}}
\nek{\rT}  {\mathbin{\sf T}}
\nek{\rtd} {\rT_2}
\nek{\rZ}  {\mathbin{\sf Z}}
\nek{\rdi} {\mathbin{{\sf D}_I}}
\nek{\rda} {\mathbin{{\sf D}_A}}
\nek{\rF}  {\mathbin{\sf F}}
\nek{\rE}  {\mathbin{\sf E}}
\nek{\rP}  {\mathbin{\sf P}}
\nek{\rS}  {\mathbin{\sf S}}

\nek{\epo} [1]
{\mathbin{{{#1}}{\vphantom{\pn}}^{\text{\bf+}}}}
\nek{\ei}   {\epo{\rE}}
\nek{\rfi}  {\epo{\rF}}

\nek{\nD}  {\mathbin{{\not\hspace{-0.35ex}\sf D}}}
\nek{\nE}  {\mathbin{{\not\hspace{-0.35ex}\sf E}}}
\nek{\nR}  {\mathbin{{\not\hspace{-0.35ex}\sf R}}}
\nek{\nF}  {\mathbin{{\not\hspace{-0.25ex}\sf F}}}
\nek{\reo} {\rE_{\hspace{-1.0pt}\Zo}}
\nek{\rez} {\rE_{\hspace{-1.0pt}\cZ}}
\nek{\nrz} {\mathbin
{\not{{\hspace{-0.4ex}\rE}}_{\hspace{-1.0pt}\cZ}}}
\nek{\nre} {\mathbin
{\not{{\hspace{-0.4ex}\rE}}_{\hspace{-1.0pt}\Zo}}}
\nek{\reff}{\rE_{\fifi}}
\nek{\reb} {\le_{\rm B}}
\nek{\eqb} {\approx_{\rm B}}
\nek{\dzo} {\dd{\Zo}}
\nek{\dde} {\dd{\rE}}
\nek{\ddec}{\dd{\Ec}}
\nek{\dep} {\dd{\rE'}}
\nek{\ddf} {\dd{\rF}}
\nek{\ddv} {\dd{\vt}}
\nek{\ddz} {\dd\cZ} 
\nek{\der} {\er-}
\nek{\dee} {\dd{\rE,\rE'}}
\nek{\Def} {\dd{\rE,\rF}}
\nek{\ddd} [2] {\dd{#1,#2}}

\nek{\fdo}{\hbox{\raisebox{0.2ex}{\mtho\tiny$\bullet$}}}
\nek{\fdt}{\hbox{\raisebox{-0.25ex}{\LARGE\bf.}}}
\nek{\bdot}[1] {\raisebox{-0.07ex}{\mtho$\stackrel{\fdt}{#1}$}}
\nek{\doa} {{\bdot a}} 
\nek{\dox} {{\bdot x}}
\nek{\dog} {\raisebox{-0.28ex}{\mtho$\stackrel{\fdt}g$}}
\nek{\doxl}{\dox_{\tt left}}
\nek{\doxr}{\dox_{\tt right}}
\nek{\rkb} [1] {|#1|_{\rm CB}}
\nek{\rkt} [2] {{^{\om}\hspace*{-1.3pt}|#1|_{#2}}}
\nek{\rko} [2] {{|#1|_{#2}}}
\nek{\rkT} [1] {{^{\om}\hspace*{-1.3pt}|#1|}}
\nek{\rkO} [1] {{|#1|}}
\nek{\kt}  [1] {{^{#1}\hspace*{-0.8pt}T}}
\nek{\bsf} [1] {I_{#1}}
\nek{\fri} [1] {\mathbin{\displaystyle\rE^{\tt fr}_{#1}}}
\nek{\fre} [1] {\frt_{{#1}}}
\nek{\fps} [3] {{\prod_{#3}{#2}\,/\,{#1}}}
\nek{\fss} [3] {{\sum_{#3}{#2}\,/\,{#1}}}
\nek{\fpt} [2] {{\sum{\sis{#2}{}}\,/\,{#1}}}
\nek{\fpd} [2] {{#1}\otimes{#2}}
\nek{\fsm} [2] {{#1}\oplus{#2}}
\nek{\dsu} {\fsm}

\nek{\sto} {[s_0,t_0]}
\nek{\ostp}{[s',t']}
\nek{\ost} {[s,t]}
\nek{\osu} {[s,u]} 
\nek{\otu} {[t,u]}
\nek{\ouv} {[u,v]}

\nek{\zs} {{\tilde s}}
\nek{\zt} {{\tilde t}}
\nek{\zu} {{\tilde u}}
\nek{\zv} {{\tilde v}}
\nek{\zn} {{\tilde n}}
\nek{\zm} {{\tilde m}}

\nek{\ff}[2] {F_{#1}^{#2}}

\nek{\spa} {\dS}
\nek{\spn} {\spa^\dN}
\nek{\spp} {\dY}

\nek{\poq} {\underline}
\nek{\nad} {\overline}
\nek{\nadd}{\widehat}

\nek{\sm} {\text{\mtho\large\boldmath$\fda$}}
\nek{\smp} [2] {\sm_{#1}^{#2-1}}
\nek{\smy} [1] {\sm_{#1}^{\iy}}

\nek{\sui} [1] {\cS_{\ans{#1}}}
\nek{\sun} {{\sui{1/n}}}
\nek{\eun} {\rS_{\ans{1/n}}}
\nek{\srn} {{\sui{r_n}}}
\nek{\ern} {\rS_{\ans{r_n}}}
\nek{\erpn} {\rS_{\ans{r'_n}}}
\nek{\suo} {{\sui0}}
\nek{\eso} {\rE_{\hspace{-1.0pt}\suo}}
\nek{\nso} {\mathbin
{\not{{\hspace{-0.4ex}\rE}}_{\hspace{-1.0pt}\suo}}}

\nek{\gal} [3] {{\tt Gal}^{#1}_{#2}(#3)} 

\nek{\nbd} [1] {{\skr O}_1(#1)}

\nek{\aH} {H^\ast}
\nek{\aB} {B^\ast}

\nek{\ovl} [1] {\overline{#1}} 
\nek{\ovg} [1] {\ovl{g_{#1}}}
\nek{\ovp} [1] {\ovl{\ga_{#1}}}

\nek{\plo} {+1} %{\!+\!1}

\nek{\yk} [1] {k_{#1}}

\nek{\mtho}{\mathsurround=0mm}
\nek{\msur}{\hspace{-1\mathsurround}}
\nek{\psur}{\hspace{0.3\mathsurround}}
\nek{\dsur}{\hspace{-0.3\mathsurround}}
\nek{\hsur}{\hspace{-0.5\mathsurround}}
\nek{\noi}{\noindent}
\nek{\vom}{\vspace{1mm}}
\nek{\vtm}{\vspace{2mm}}

\nek{\uv}{{\bf V}}

\nek{\wA} {{\widehat A}}
\nek{\ha} {{\hat a}}
\nek{\hb} {{\hat b}}
\nek{\he} {{\hat\ve}}
\nek{\hT} {{\hat t}}
\nek{\hl} {{\hat l}}
\nek{\hs} {{\hat s}}
\nek{\hx} {{\hat x}}
\nek{\hm} {{\hat m}}
\nek{\hn} {{\widehat n}}
\nek{\hk} {{\widehat k}}

\nek{\fo} {{\mathbf 0}}
\nek{\fr} {{\mathbf 1}}

\nek{\vnu} {{\vec \nu}} 
\nek{\dvn} {\dd\vnu} 
\nek{\wnu} {\cW_{\vnu}} 
\nek{\ewn} {\rE_{\vnu}}
\nek{\ret} {\rE_{T}}

\nek{\ci} [1] {I_{#1}}

\nek{\vex}{{\vec x}}
\nek{\vey}{{\vec y}}

\nek{\dns} [1] {d_{\vnu}^{#1}}
\nek{\den} {d_{\vnu}}
\nek{\din} {d_\vnu}

\nek{\poo}{=_{\tt df}}

\nek{\dpi}{d_\vpi}

\nek{\nol}[1] {{\bf 0}_{#1}}
\nek{\edi}[1] {{\bf 1}_{#1}}

\nek{\nrn} {_{\ans{r_n}}}
\nek{\vpr} {\vpi\nrn}
\nek{\dpr} {d\nrn}
\nek{\drn} {\dd{\ans{r_n}}}

\nek{\okr} [2] {{\skr O}_{#1}(#2)}

\nek{\nid}{\gE}

\nek{\zid} [2] {\cI_{#1/#2}}
\nek{\zidef} {\zid\rE\rF}

\nek{\zfo} [2] {\dP_{#1/#2}}
\nek{\zfoef} {\zfo\rE\rF}

\nek{\pei} {\dP_{\Eo\Ei}}
\nek{\pet} {\dP_{\Eo\Et}}
\nek{\peo} {\dP_{\Eo}}
\nek{\pep} {\dP'_{\Eo}}
\nek{\iet} {\cI_{\Eo\Et}}
\nek{\iei} {\cI_{\Eo\Ei}}
\nek{\ieo} {\cI_{\Eo}}

\SetMathAlphabet{\cur}{bold}{U}{eur}{b}{n}
\renek{\Gd} {\fG_\fda}
\renek{\Fs} {\fF_\fsg}

\nek{\ubf}{\fontseries{b}\selectfont}

\mathchardef\alphA ="710B
\mathchardef\betA ="710C
\mathchardef\gammA ="710D
\mathchardef\deltA ="710E
\mathchardef\vartA ="7123
\mathchardef\kpA   ="7114
\mathchardef\mU    ="7116
\mathchardef\nU    ="7117
\mathchardef\pI    ="7119
\mathchardef\rhO   ="711A
\mathchardef\sigmA ="711B
\mathchardef\taU   ="711C
\mathchardef\xI ="7118

\nek{\fal} {{\cur{\alphA}}}
\nek{\fba} {{\cur{\betA}}}
\nek{\fsg} {{\cur{\sigmA}}}
\nek{\fpi} {{\cub{\pI}}}
\nek{\fta} {{\cub{\taU}}}
\nek{\fda} {{\cur{\deltA}}}

\nek{\dds}{\dd\fsg}

\nek{\nmp} {\Longleftarrow}

\nek{\tsc}[1]{\hbox{\footnotesize\sc{#1}}}

\nek{\ddi} {\dd{\cI}}
\nek{\ddj} {\dd{\cJ}}
\nek{\ddij}{\dd{(\cI,\cJ)}}

\nek{\ren}  {\le_{\tsc{c}}}
\nek{\rens} {<_{\tsc{c}}}
\nek{\eqc}  {\sim_{\tsc{c}}}

\nek{\eui}[1] {{\text{\it{EU}}}_{\ans{#1}}}

\nek{\sqa} {\sq^\ast}
\nek{\sua} {\su^\ast}

\nek{\prf} [1] {\hbox{\S\hspace{0.3ex}\ref{#1}}}
\nek{\prff}[1] {\S\S~\ref{#1}}
\nek{\pff} [1] {\S~{#1}}

\nek{\mrn} {\mu\nrn}

\renek{\dop} [1] {\complement #1}
\renek{\df}  [1] {{#1}^\complement}
\nek{\doP}  [1] {{#1}^\complement}

\nek{\ima} [2] {#1\text{\hspace*{0.1ex}''}#2}
\nek{\imb} [2] {#1\text{\hspace*{0.3ex}''\hspace*{-0.2ex}}#2}

\nek{\aprb}{\approx_{\tsc b}}
\nek{\ismb}{\cong_{\tsc b}}

\nek{\incs} {<_{\tsc i}}
\nek{\inc} {\le_{\tsc i}}
\nek{\eqi} {\sim_{\tsc i}}

\nek{\reas} {<_{\tsc a}}
\nek{\rea} {\le_{\tsc a}}
\nek{\eqa} {\sim_{\tsc a}}

\nek{\reaas} {<_{\tsc{aa}}}
\nek{\reaa} {\le_{\tsc{aa}}}
\nek{\eqaa} {\sim_{\tsc{aa}}}

\nek{\rebs} {<_{\tsc b}}
\renek{\reb} {\le_{\tsc b}}
\renek{\eqb} {\sim_{\tsc b}}

\nek{\reBs}{<_{\tsc{bm}}}
\nek{\reB} {\le_{\tsc{bm}}}
\nek{\eqB} {\sim_{\tsc{bm}}}

\nek{\rds} {<^{\Da}_{\tsc{b}}}
\nek{\rd} {\le^{\Da}_{\tsc{b}}}
\nek{\eqd} {\approx^{\Da}_{\tsc{b}}}

\nek{\rbas} {<_{\tsc{b,ba}}}
\nek{\rba} {\le_{\tsc{b,ba}}}
\nek{\eqba} {\approx_{\tsc{b,ba}}}

\nek{\rbasp} {<_{\tsc{b,ba}}^+}
\nek{\rbap} {\le_{\tsc{b,ba}}^+}
\nek{\eqbab} {\approx_{\tsc{b,ba}}^+}

\nek{\nab} [1] {\nabla(#1)}
\nek{\orb} {\le_{\tsc{rb}}}
\nek{\srb} {<_{\tsc{rb}}}
\nek{\erb} {\sim_{\tsc{rb}}}

\nek{\orbpp} {\le_{\tsc{rb}}^{++}}
\nek{\srbpp} {<_{\tsc{rb}}^{++}}
\nek{\erbpp} {\sim_{\tsc{rb}}^{++}}

\nek{\orbp} {\le_{\tsc{rb}}^{+}}
\nek{\srbp} {<_{\tsc{rb}}^{+}}
\nek{\erbp} {\sim_{\tsc{rb}}^{+}}

\nek{\ork} {\le_{\tsc{rk}}}
\nek{\srk} {<_{\tsc{rk}}}
\nek{\erk} {\sim_{\tsc{rk}}}

\nek{\obe} {\le_{\tsc{be}}}
\nek{\sbe} {<_{\tsc{be}}}
\nek{\ebe} {\sim_{\tsc{be}}}

\nek{\rei} {\rE_{\cI}}
\nek{\rej} {\rE_{\cJ}}

\nek{\eeb} {\sim_{\tsc{b}}}

\nek{\obep} {\le_{\tsc{be}}^+}
\nek{\sbep} {<_{\tsc{be}}^+}
\nek{\ebep} {\sim_{\tsc{be}}^+}

\nek{\seb} {<_{\tsc{b}}}

\nek{\supp} {\mathop{\tt supp}}

\nek{\atp} {\mathop{\tt at}^+}
\nek{\atm} {\mathop{\tt at}^-}

\nek{\hv} [2] {||#1||_{#2}}

\nek{\vmu} {{\vec \mu}}

\nek{\bel} [1] {\mathrel{\text{\boldmath\mtho$\ell$}}^{#1}}
\nek{\beL} [1] {\mathrel{\text{\bfsf L}}^{#1}}
\nek{\bem}     {\mathrel{\text{\bfsf m}}}
\nek{\fCo}{\mathbin{\rC_{\text{\sf0}}}}
\nek{\fco}{\mathrel{\text{\bfsf c}_{\text{\sf0}}}}
\nek{\fc} {\mathrel{\text{\bfsf c}}}
\nek{\fvt}{\mathrel{\text{\boldmath\mtho$\vt$}}}

\nek{\beli} {\bel\iy}

\nek{\oin} {{[0,1]^\dN}}

\nek{\iz} [2] {\cI_{#1}^{#2}}
\nek{\jz} [1] {\iz{}{}(#1)}
\nek{\iw} [2] {\cW_{#1}^{#2}}
\nek{\jw} [1] {\iw{}{}(#1)}
\nek{\ib} [2] {\cB^{#1}_{#2}}
\nek{\jb} [1] {\ib{}{#1}}

\nek{\ovu} {{\overline u}}

\nek{\nr}[2]{\nor{#1}_{#2}} %{||#1||_{#2}}

\nek{\TS}{\textstyle}
\nek{\DS}{\displaystyle}

\nek{\Ba}{B^\ast}

\nek{\resi} [1] {\mathop{\restriction_{#1}}}

\nek{\gP}{{\BBB{\mathfrak P}\BBB}}
\nek{\gF}{{\BBB{\mathfrak F}\BBB}}
\nek{\gJ}{{\BBB{\mathfrak J}\BBB}}
\nek{\dG}{{\BBB{\mathbb G}\BBB}}
\nek{\dH}{{\BBB{\mathbb H}\BBB}}

\nek{\ac} {\cdot} %action

\nek{\curle}{\preccurlyeq}
\nek{\cle}{\curle}
\nek{\cge}{\succcurlyeq}
\nek{\cl} {\prec}
\nek{\resic} [1] {\resi{\cle #1}}
\nek{\rec} {\resic}
\nek{\rc} {\rsic}
\nek{\recb} [1] {\rec{(#1)}}
\nek{\rsic} [1] {\resi{\cl #1}}
\nek{\rcb} [1] {\rc{(#1)}}

\nek{\kai} {\forall^\iy\hspace{0.1ex}}
\nek{\exi} {\exists^\iy\hspace{0.1ex}}

\nek{\ovw}{\nad w}

\nek{\il}[2] {\ir{n_{#1}}{n_{#2}}}
\nek{\ia}[2] {\ir{a_{#1}}{a_{#2}}}
\nek{\ij}[2] {\ir{j_{#1}}{j_{#2}}}
\nek{\cB}{{\BBB{\skr B}\BBB}}
\nek{\cG}{{\BBB{\skr G}\BBB}}
\nek{\cL}{{\BBB{\skr L}\BBB}}
\nek{\cN}{{\BBB{\skr N}\BBB}}
\nek{\cU}{{\BBB{\skr U}\BBB}}

\nek{\pnd}{\pn^\dN}
\nek{\dnd}{{(\dn)}{}^\dN}

\nek{\anp} [2] {\ang{#1}^{#2}}

\nek{\ta}{\tau}

\nek{\lev} {\mathop{\tt{lev}}}
\nek{\glu} {\mathop{\tt{dep}}}
\nek{\dia} {\mathop{\tt{diam\hspace{0.15ex}}}}

\nek{\Ei}{\rE_{\text{\sf 1}}}
\nek{\Ed}{\rE_{\text{\sf 2}}}
\nek{\Et}{\rE_{\text{\sf 3}}}
\nek{\Ey}{\rE_\iy}
\nek{\Eya} {\mathbin{(\rE_\iy)^\alo}}

\nek{\npi}{\nu_\vpi}
\nek{\nsi}{\nu_\psi}

\nek{\Ii}{\cI_1}
\nek{\Id}{\cI_2}
\nek{\It}{\cI_3}

\nek{\emb}  {\sqsubseteq_{\tsc{b}}}
\nek{\emn}  {\sqsubseteq_{\tsc{c}}}
\nek{\embi} {\sqsubseteq_{\tsc{b}}^{\rm i}}
\nek{\emni} {\sqsubseteq_{\tsc{c}}^{\rm i}}

\nek{\sio}{\cS_0}
\nek{\esn}{\rS_{\ans{1/n}}}

\nek{\nsn}{\mathbin{\not\hspace*{-0.3ex}\esn}}

\nek{\req}[2]{{\DS|_{#1}^{#2}}}

\nek{\rlq}[1]{\resi{\ge#1}}
\nek{\rmq}[1]{\resi{<#1}}
\nek{\rnq}[1]{\resi{\le#1}}

\nek{\Dij} {\dd{\cI,\cJ}}
\nek{\deo} {\dd{\Eo}}

\nek{\dc}[2] {{\bf W}^{#1}_{#2}}

\nek{\Za}{{\nad Q}}
\nek{\Ya}{{\widetilde Y}}

\nek{\Ga}{\Gamma}
\nek{\rG}  {\mathbin{\sf G}}

\nek{\inva}{{\tt{inv}}}

\nek{\tP}{{P^\ast}}
\nek{\app} {{\hspace{0.2ex}{\cdot}\hspace{0.2ex}}}

\nek{\lland}{\,\land\,}

\nek{\seqv} {\hspace{0.3ex}\Leftrightarrow\hspace{0.3ex}}
\nek{\simp} {\hspace{0.3ex}\Rightarrow\hspace{0.3ex}}

\nek{\eqn} [1] {\equiv_{#1}}
\nek{\eqk} [2] {\equiv_{#1}^{#2}}
\nek{\eqo} [2] {\xrightarrow{#1,#2}}

\nek{\pit}{\tilde\pi}
\nek{\tve}{\tilde\ve}

\nek{\ef} {\dd{\rE,\rF}}

\nek{\isi} {\cong}
\nek{\isa} {\mathrel{\hspace{0.2ex}\isi^\ast}}

\nek{\ske} [3] {\equiv_{#1#2}^{#3}}
\nek{\skab}[1] {\ske AB{#1}}
\nek{\spab}[1] {\ske{A'}{B'}{#1}}

\nek{\Indent}{\hspace*{3ex}}
\nek{\fsur}{\hspace{0.5\mathsurround}}

\nek{\rdm}  {{\mathsf c}_{\tt max}}

\nek{\co} {\dd{\text{\bfsf c}_{\mathsf0}}}
\nek{\lv} {{\normalfont\scshape lv}-}

\nek{\susi} {\exists^{\iy}\hspace*{0.2ex}}
\nek{\kazi} {\forall^{\iy}\hspace*{0.2ex}}

\nek{\igi}{$\hbox{\mtho\boldmath$\fal$}$}
\nek{\igd}{$\hbox{\mtho\boldmath$\fba$}$}

\nek{\cg}[1] {{\tt Choq}(#1)}
\nek{\cgs}[1]{{\tt Choq}^{\rm s}(#1)}

\nek{\dnn}{\dN^\dN}
\nek{\nnn}{(\dnn){\vphantom{\dN}}^\dN}

\nek{\isg} {{S_\iy}}

\nek{\uset}{{Universal sets}}

\nek{\ler}[2] {\mathbin{\sim^{#1}_{#2}}}
\nek{\aer}[2] {\mathbin{\rE_{#1}^{#2}}}
\nek{\ergx}{\aer\dG\dX}
\nek{\egx} {\ergx}
\nek{\lo} [3] {\cO(#3,#1,#2)}
\nek{\sym}{\ler} % [2] {\sim_{#1}^{#2}}
\nek{\ong} {1_\dG}
\nek{\rr} [2] {\rR^{#1}_{#2}}

\nek{\rav}[1]
{\mathbin{\mathsf\Delta}_{#1}}
\nek{\toq}{\circle*{0.5}}
\nek{\tob}{\circle*{1.0}}

\nek{\stob}{\circle{3.0}}
\nek{\ktob}{\kras\circle{3.0}}

\nek{\cob}{\circle{1.5}}
\nek{\mtir}{\line(-1,0){2}}

\nek{\bon} [2] {\cO_{#1}(#2)}

\nek{\dnnn} {\dnnp\dN}
\nek{\dnnp} [1] {(\dnn){\vphantom{\dN}}^{#1}}

\nek{\prift}{\sf}
\nek{\Penu}{{\prift$\id11$ Enumeration}}
\nek{\Refl}{{\prift Reflection}}
\nek{\Uset}{{\prift Universal Sets}}
\nek{\Kres}{{\prift Kreisel Selection}}
\nek{\Cenu}{{\prift Countable-to-1 Enumeration}}
\nek{\Cuni}{{\prift Countable-to-1 Uniformization}}
\nek{\Cpro}{{\prift Countable-to-1 Projection}}
\nek{\Bore}{{\prift Borel Extension}}
\nek{\Sepa}{{\prift Separation}}
\nek{\Redu}{{\prift Reduction}}

\nek{\dpf} {\mathord{{\dP^2\hspace*{-0.3ex}}\res\rF}}
\nek{\dpe} {\mathord{{\dP^2\hspace*{-0.3ex}}\res\rE}}

\nek{\ek}[2] {[#1]_{{#2}}}

\nek{\eke}[1] {\ek{#1}{\rE}}
\nek{\ekeo}[1] {\ek{#1}{\Eo}}
\nek{\ekec}[1] {\ek{#1}{\Eco}}
\nek{\ekco}[1] {\ek{#1}{\Ec}}
\nek{\ekfo}[1] {\ek{#1}{\Fo}}
\nek{\ekf}[1] {\ek{#1}{\rF}}
\nek{\ekg}[1] {\ek{#1}{G}}

\nek{\ur}{_{\tt right}}
\nek{\ul}{_{\tt left}}

\nek{\cont}{\hbox{\mtho\large$\gc$}}

\nek{\mem} {\dd\in}

\nek{\bV}{{\mathbf V}}
\nek{\bL}{{\bf L}}
\nek{\cli} {{\sc cli}}

\nek{\di} [1] {{#1}^\#}
\nek{\drE} {\mathbin{\di\rE}}
\nek{\drF} {\mathbin{\di\rF}}

\nek{\kon} {\hbox{\mtho\large${\mathfrak c}$}}
\nek{\bk} [1] {{\cur B}_{#1}}

\nek{\wtau}{{\widehat\tau}}

\nek{\gra}[1]{\dd{#1}grainy}
\nek{\grap}{grainy}
\nek{\Grap}{Grainy}

\nek{\xE}[2] {\mathbin{\rR^{#2}_{\ge #1}}}

\nek{\moq} [1] {\Mod_{#1}}
\nek{\mox} {\moq} %[1] {{\mathbb X}_{#1}}
\nek{\loa} [1] {j_{#1}}
\nek{\ism} [1] {\cong_{#1}}
\nek{\izm} [2] {\cong_{#1}^{#2}}
\nek{\aut} [1] {\Aut_{#1}}

\nek{\hfn} {{\rm HF}(\dN)}
\nek{\tce} [1] {{\rm TC}_\ve(#1)}
\nek{\ihf} {\simeq_{\hfn}}

\nek{\symr}{\equiv}
\nek{\rrt} [3] {\symr_{#2#3}^{#1}}
\nek{\rrq} [5] {{#4}\symr_{#2#3}^{#1}{#5}}
\nek{\nrq} [5] {{#4}\not\symr_{#2#3}^{#1}{#5}}
\nek{\rrQ} [5] {{#4}\symr_{#2\,,\,#3}^{#1}{#5}}
\nek{\rro} [5] {\ang{#2,#4}\symr^{#1}\ang{#3,#5}}
\nek{\rrO} [1] {\symr^{#1}}

\nek{\lww} {\cL_{\omi\om}}

\nek{\forc}
{\mathrel{{|\hspace{-0.2ex}|\hspace{-1.1ex}-\hspace{-1.7ex}-}}}

\newlength{\dxii}
\nek{\fC} {{\bf C}}
\nek{\pg} {\fC_\dG}
\nek{\px} {\fC_\dX}

\nek{\gen} {gen.\ }
\nek{\hg}  {loc.\hspace{0.4ex}gen.\ }
\nek{\hgp} {loc.\hspace{0.4ex}gen.}
\nek{\rfy} {\rF^\iy}

\nek{\incl} [1] {\mathop{\text{\sc Int}}\overline{#1}}  %{{\tt IntCl}}

\nek{\PP} {pinned}
\nek{\PPP}{Pinned}

\nek{\bap}{{\bar p}}
\renek{\wtau} {{\widehat p}}

\nek{\lbr} {\linebreak[0]}

\nek{\zO} {\yo}
\nek{\zi} {\yi}
\nek{\zd} {\yd}
\nek{\zT} {\yt}
\nek{\yo} {,\linebreak[0]}
\nek{\yi} {,\linebreak[0]\hspace{0.1ex}}
\nek{\yd} {,\linebreak[0]\:}
\nek{\yt} {,\linebreak[0]\;}
\nek{\zq} {,\linebreak[0]\;\,}

\nek{\prit} [1] {[{{\rm #1}}]}
\nek{\nor} [1] {\|#1\|}

\nek{\fap} {f.\hspace{0.1ex}a.\hspace{0.1ex}p.\hspace{0.1ex}m.}

\nek{\bpr} [1] {\bpf[{{\sl#1}\/}]}

\nek{\eqr} {equivalence relation}

\nek{\fx} {{\bf x}}

\nek{\rzd} {^{\text{\tt red}}}

\nek{\srez} [2] {{\bf S}_{#1}(#2)}

\nek{\qc}  [1] {\resi{> #1}}
\nek{\qec} [1] {\resi{\ge #1}}
\nek{\rme} [1] {\resi{\le#1}}

\nek{\alex} {<_{\text{\tt alex}}}
\nek{\lex} {<_{\text{\tt lex}}}
\nek{\act} {<_{\text{\tt act}}}

\nek{\bbo} {\mathbb 0}

\nek{\fz} {{\mathbf z}}

\nek{\dln} {2\lom}

\nek{\fK} {{\bf K}}

\nek{\Ks} {\fK_\fsg}

\nek{\dva}{{\ans{0,1}}}

\nek{\rtn}{\dR^\dN}
\nek{\ntn}{\dN^\dN}
\nek{\ztn}{\dZ^\dN}

\nek{\snos} [1] {\,\footnote{\ #1}}
\nek{\snom}   {\,\footnotemark}
\nek{\snot} [1] {\footnotetext{\ #1}}

\nek{\rH}  {\mathbin{\sf H}}

\nek{\fras}[2] {\text{\footnotesize$\DS\frac{#1}{#2}$}}
\nek{\fral}[2] {\text{\large$\frac{#1}{#2}$}}

\nek{\renu}{\tenu{{\rm(\roman{enumi})}}}

\nek{\ergg}{\aer\dG\dG}

\nek{\rit} [1] {{\it#1\/}}

\nek{\lap} [1] {``#1''}

\nek{\mto} {\mapsto}

\nek{\pgcr} {\addtocounter{page}1}

\nek{\gs} [2] {g_{#1#2}}

\nek{\bie} [2]
{\mathord{{\text{\large\mtho$\prec$}}#1,#2{\text{\large\mtho$\succ$}}}}
\nek{\wa} {\widehat a}
\nek{\wb} {\widehat b}
\nek{\ws} {\bar s}
\nek{\wt} {\bar t}

\nek{\itlm}[1] {\itla{#1}\hspace*{0.0ex}\imar{#1}}

\nek{\xH} {\mathbf F}
\nek{\xh} [1] {\xH_{>#1}}
\nek{\yh} [1] {\xH_{\le#1}}
\nek{\zh} [1] {\xH_{\ge#1}}

\nek{\nksp} {\hspace{0.1ex}}
\nek{\hh} [2] {{\DS H_{#1\nksp#2}}}
\nek{\jf} [2] {{\DS {\mathbf A}^{\mathtt{fin}}_{#1\nksp#2}}}
\nek{\af} [2] {{\DS \cA^{\mathtt{fin}}_{#1\nksp#2}}}
\nek{\Af} [2] {{\DS U^{\mathtt{fin}}_{#1\nksp#2}}}
\nek{\ji} [2] {{\DS {\mathbf A}^{\mathtt{1}}_{#1\nksp#2}}}
\nek{\ai} [2] {{\DS \cA^{\mathtt{1}}_{#1\nksp#2}}}
\nek{\Ai} [2] {{\DS U^{\mathtt{1}}_{#1\nksp#2}}}
\nek{\ds} {{\mathbf S}}
\nek{\dq} {{\mathbf H}}
\nek{\gru} {\hspace*{0.2ex}}

\nek{\lmt} [1] {\le_{\text{\sc nt}}^{#1}}
\nek{\emt} [1] {\rE^{\vphantom{\xi}#1}_{\text{\sc nt}}}
\nek{\mti} [1] {\cJ_{\text{\sc nt}}^{#1}}
\nek{\wmti}[1] {{\wcI}_{\text{\sc nt}}^{#1}}

\nek{\scw} {<_{\text{\tt cw}}}
\nek{\gcw} {\ge_{\text{\tt cw}}}
\nek{\lcw} {\le_{\text{\tt cw}}}
\nek{\lnt} {\le_{\text{\sc nt}}}
\nek{\ent} {\rE_{\text{\sc nt}}}
\nek{\nti} {\cJ_{\text{\sc nt}}}
\nek{\wnti}{\wcJ_{\text{\sc nt}}}

\nek{\ndi}  {\cI_{\text{\sc nt}}}
\nek{\ndie} {\cE_{\ndi}}

\nek{\bez} {\dif}

\nek{\Eqr} {Equivalence relation}

\nek{\rah} [1] {\text{\tt rnk}({#1})}
\nek{\raf} [1] {\text{\tt rnk}_{#1}}
\nek{\rag} [2] {\text{\tt rnk}_{#1}(#2)}

\nek{\pcw} {\mathbin{+_{\text{\hspace*{-0.5ex}\tt cw}}}}

\nek{\tlw} {2\lom}
\nek{\nw} {\dN^\om}
\nek{\nlw} {\dN\lom}
\nek{\tw} {2^\om}

\nek{\nt} {\text{\ubf{NT}}}

\nek{\wS} {\widehat S}
\nek{\wQ} {\widehat Q}

\nek{\ceq} [3] {{#1}_{#2#3}}
\nek{\cer} [4] {{#1}_{#2#3}^{#4}}

\nek{\rU}  {\mathbin{\sf U}}

\nek{\bai} {\cN}

\nek{\lip} [2] {\text{\sc Emb}(#1,#2)}
\nek{\wid} [1] {\text{\sc Wid}(#1)}
\nek{\vid} [1] {\text{\sc Wid}(#1)}

\nek{\lmto}{\longmapsto}

\nek{\atlh} {\addtolength{\itemsep}{-\dxh}}
\nek{\atli} {\addtolength{\itemsep}{-\dxi}}

\nek{\cT} {{\skr T}}

\newlength{\dxh}
\nek{\cip} {{\sc cip}} %{c.\,i.\,p.}

\nek{\uap} [2] {A^{#1}(#2)}
\nek{\uaq} [3] {A^{#1}_{#3}(#2)}
\nek{\use} [1] {A(#1)}

\theoremstyle{definition}
\newtheorem{coj}       [theore]{Conjecture} 
\nek{\bcj} {\begin{coj}}
\nek{\ecj} {\end{coj}}

\nek{\api}{\addtocounter{page}{1}}

\nek{\cntr} [1] {{#1}^{\downarrow}}

\nek{\gam} [3] {\text{\ubf G}_{#1}^{#2}(#3)}

\nek{\ts} {\widehat s}
\nek{\tsg}{\widehat \sg}

\nek{\imae}[1]{}

\vyk{
\nek{\imae}[1]{\marginpar[%\vspace{-1ex}%
\flushright\footnotesize\vspace{-4ex}%
$\mtho\longrightarrow$\\%
\vspace{-1ex}{#1}$\mtho$\vspace*{2ex}]%
{%\vspace{-1ex}%
\flushleft\footnotesize\vspace{-4ex}%
$\mtho\longleftarrow$\\%
\vspace{-1ex}{#1}$\mtho$\vspace*{2ex}}
}%
}

\nek{\lam} [1] {\label{#1}\hspace*{-3pt}\imar{#1}}%
\nek{\las} [1] {\label{#1}\imae{#1}}%
\nek{\urav}  [1] {\bus\itsep#1\eus} 

\nek{\efp}{\dd{\rE,{\rF'}}}

\nek{\dfp} {\dd{\rF'}}

\nek{\inw} [2] {\dd{(#1\to#2)}invariant}

\nek{\come} {co-meager\ }

\nek{\bfei} {{\ubf either}}
\nek{\bfor} {{\ubf or}}

\nek{\gla} {Chapter}
\nek{\paf} {Section}

\nek{\grf} [1] {\gla~\ref{#1}}
\nek{\nrf} [1] {\paf~\hbox{\ref{#1}}}

\nek{\dPP} {\dP\ti\dP}
\nek{\dSP} {\dS\ti\dP}
\nek{\xib} {{\bma\xI\hspace{0.1ex}}}
\nek{\sgb} {\bma\sigmA}
\nek{\ddp} {\dd\dP}
\nek{\coll} {\text{\sc Coll}}
\nek{\tal}{\xib\ul}
\nek{\tar}{\xib\ur}
\nek{\dCP} {\dC\ti\dP}
\nek{\dCC} {\dC\ti\dC}
\nek{\sgbl}{\sgb\ul}
\nek{\sgbr}{\sgb\ur}

\nek{\uoa} {\mathbin{{\text{\bfsf u}}^\ast_0}}
\nek{\uo}  {\mathbin{{\text{\bfsf u}}_0}}

\nek{\ld} [1] {\mathbin{{\rD}(#1)}}

\nek{\xd} [2] {{\mathrel{{\td}(#1\hspace*{0.1ex};\hspace*{0.1ex}#2)}}}
\nek{\qd} [3] {{\mathrel{{\td}
(\ang{#1\hspace*{0.1ex};\hspace*{0.1ex}#2}_{#3\in\dN})}}}
\nek{\qqd} [3] {{\mathrel{{\td}
(\ang{#1\hspace*{0.1ex};\hspace*{0.1ex}#2}_{#3})}}}

\nek{\td}  {\mathrel{\rD}}
\nek{\ntd} {\mathrel{{\not\hspace{-0.0ex}\td}}}
\nek{\mast} {\text{\mtho$\ast$}}
\nek{\ja} {\text{\mtho\boldmath$a$}}
\nek{\jx}{\mathbf{x}}
\nek{\cemu}{

\itsep
}

\nek{\aemu}{

\itsep
}

\nek{\xs}[2] {\cF_{#1}^{#2}}
\nek{\xb}[2] {\mathord{^{\text{\tt T}}\hspace*{-0.7ex}\xs{#1}{#2}}}
\nek{\rnd} {\rn\hspace{-0.4ex}\ti\hspace{-0.3ex}\rn}
\nek{\rnn} [1] {\drb{#1}\hspace{-0.4ex}\ti\hspace{-0.3ex}\rn}
\nek{\rnne}[1] {\drbe{#1}\hspace{-0.4ex}\ti\hspace{-0.3ex}\rn}
\nek{\Eit} {\rE_{13}}
\nek{\Ejt} {\Eit}

\nek{\rg}  [1] {\resi{> #1}}
\nek{\rl}  [1] {\resi{< #1}}
\nek{\rge} [1] {\resi{\be #1}}
\nek{\rle} [1] {\resi{\me #1}}

\nek{\resa} [1] 
{\mathop{{\mtho\text{\large$\restriction$}}}\DS_{#1}}
\nek{\reza} [1] 
{\mathop{{\mtho\text{\large$\restriction$}}}\DS^{#1}}
\nek{\resb} [2] 
{\mathop{{\mtho\text{\Large$\restriction$}}}\DS_{#1}^{#2}}
\renek{\resi} [1] 
{\mathop{\TS\hspace{-0.3ex}\restriction\hspace{-0.2ex}}\DS_{#1}}

\nek{\rga}  [1] {\resa{> #1}}
\nek{\rgea} [1] {\resa{\be #1}}

\nek{\xgg}  [2] {\resb{> #1}{>#2}}
\nek{\xge}  [2] {\resb{> #1}{\be#2}}
\nek{\xeg}  [2] {\resb{\be #1}{>#2}}
\nek{\xee}  [2] {\resb{\be #1}{\be #2}}

\nek{\rre} [1] {\dR^{\be#1}}
\nek{\rrs} [1] {\dR^{>#1}}
\nek{\rse} [1] {\dR^{\me#1}}
\nek{\rss} [1] {\dR^{<#1}}

\nek{\drm} [1] {\dR^{<#1}}
\nek{\drb} [1] {\dR^{>#1}}
\nek{\drme} [1] {\dR^{\me#1}}
\nek{\drbe} [1] {\dR^{\be#1}}

\nek{\xgo}  [1] {\resa{> #1}}
\nek{\xeo}  [1] {\resa{\be #1}}
\nek{\xog}  [1] {\reza{> #1}}
\nek{\xoe}  [1] {\reza{\be #1}}

\nek{\zgo}  [1] {\resa{<#1}}
\nek{\zeo}  [1] {\resa{\me #1}}  

\nek{\be}{\geqslant}
\nek{\me}{\leqslant}

\nek{\kom} [2] {(#1)_{#2}}
\nek{\kmi} [1] {\kom{#1}1}
\nek{\kmd} [1] {\kom{#1}2}

\nek{\ekit} [1] {\ek{#1}{{\Eit}}}
\nek{\ekt} [1] {\ek{#1}{{\Et}}}
\nek{\etk} [2] {\ek{#1}{{\eqn{#2}}}}

\renek{\ek}[2] {[\hspace*{0.2ex}#1\hspace*{0.2ex}]_{{#2}}}

\nek{\ig} [2] {#1(#2)}
\nek{\ih} [2] {#1[#2]}

\nek{\uh} [2] {#1(#2)}
\nek{\wh} [2] {#1(\ekt{#2})}

\nek{\jh} [4] {#1_{#2}^{#4}[#3]}
\nek{\ah} [4] {\widehat{#1}_{#2}^{#4}[#3]}

\nek{\jg} [4] {#1_{#2}^{#4}(#3)}
\nek{\ag} [4] {\widehat{#1}_{#2}^{#4}(#3)}

\nek{\vh} {\jh} %[4] {#1_{#2}[#4]({\ekt #3})}
\nek{\vg} {\jg} %[4] {#1_{#2}[#4](#3)}
\nek{\wE} {\widehat E}

\nek{\zf} [3] {S_{#2}^{\ekt{#3}}(#1)}
\nek{\zfi}[2] {S_{#1}^{\ekt{#2}}}

\nek{\ov}{\overline}

\nek{\qz} [3] {\pi_{#1}^{#2}[#3]}
\nek{\ez} [2] {\ve^{#1}[#2]}
\nek{\kz} [2] {\ka^{#1}[#2]}

\renek{\pz} [2] {p^{#1}[#2]}
\nek{\zp} [3] {p_{#2}[#3](#1)}

\nek{\ze} [3] {{\ve^{#2}_{#1}[#3]}}
\nek{\zk} [3] {{\ka^{#2}_{#1}[#3]}}

\renek{\zs} [3] {{\DS s^{#2}_{#1}(#3)}}

\nek{\zka} [4] {{k_{#1}^{#4,#2}[#3]}}
\nek{\zep} [4] {{e_{#1}^{#4,#2}[#3]}}

\nek{\dk} [2] {W^{#2}_{#1}}
\nek{\dl} [1] {W^{#1}}

\nek{\eo} [1] {\mathbin{\DS{\Eo}^{#1}}}

\nek{\eRe} {\end{rem}} 

\nek{\ifd}  {{\cur{Fin}_2}}

\nek{\fje}[1] {\ifi_{\me #1}}
\nek{\fj} [1] {\ifi_{<#1}}

\theoremstyle{definition}
\newtheorem{bla}      [theore]{Blanket Agreement}
\nek{\bbl}{\begin{bla}}
\nek{\ebl}{\qed\end{bla}}
\newtheorem{conj}      [theore]{Conjecture}
\nek{\bgi}{\begin{conj}}
\nek{\egi}{\qed\end{conj}}

\nek{\fSg}{{\mathbf\Sg}}
\nek{\fPi}{{\mathbf\Pi}}

\nek{\xet} {\dd{\xi\text{-}\Et}}
\nek{\xeq} [1] {\dd{\xi\text{-}{\eqn #1}}}

\nek{\pp} [3] {{\DS\vpi_{#1}^{#2}[#3]}}

\nek{\ale} [1] {\alex^{#1}}

\nek{\jc} [2] {C^{#1}[#2]}
\nek{\jd} [2] {D^{#1}[#2]}
\nek{\jdp} [2] {D^{#1}[#2]^+}
\nek{\jp} [2] {\Pi^{#1}[#2]}

\nek{\und} {{\boldsymbol\bot}}

\nek{\kru} [3] {f^{#1\hspace{0.1ex}#2}(#3)}
\nek{\krv} [2] {f^{#1\hspace{0.1ex}#2}}

\nek{\zg} [3] {G^{#2}_{#1}[#3]}
\nek{\yg} [2] {G^{#1}[#2]}
\nek{\yy} [2] {G^{#1}_{\text{\tt fin}}[#2]}

\nek{\enp} [2] 
{\mathrel{\hspace*{0.3ex}\cong_{#1}^{#2}
\hspace*{0.1ex}}}

\nek{\enq} [3] {\mathrel{\enp{#1}{{#2\res#3}}}}

\nek{\pro} [2] {\pi^{#2}_{#1}}
\nek{\Pro} [2] {\Pi^{#2}_{#1}}

\nek{\ikl} [2] {C^{#2}_{#1}}

\nek{\ppu} {\subsubsection}

\nek{\xg} [2] {F^{#1}_{#2}}

\nek{\hhh} [2] {h^{#2}_{#1}}

\nek{\pip} [4] {\fta_{#1}^{#2}(#3,#4)}
\nek{\piq} [2] {\fta_{#1}^{#2}}

\nek{\kk} [1]  {\ka_{#1}}

\nek{\xda} [2] {\fda_{#1#2}}

\nek{\rnD} {(\rn){}^D}

\nek{\gr} [6] {G_{#1#2}^{#3#4}(#5,#6)}
\nek{\gq} [5] {G_{#1#2}^{#3#4}(#5)}
\nek{\gp} [4] {G_{#1#2}^{#3#4}}

\nek{\hq} [3] {[#3]_{#1#2}}
\nek{\hp} [2] {g_{#1#2}}

\nek{\bPi} {\mathbf\Pi}

\nek{\tiw} {{\tilde w}}
\nek{\tq} {{\tilde q}}

\nek{\po}{P''_0}

\nek{\dsf} {\dS_{\text{\tt fin}}}
\nek{\rsf} {\rS_{\text{\tt fin}}}

\nek{\bdo} {\mathbf\Pi}

\vyk{
\lccode`\-=`\-
\defaulthyphenchar=127
\usepackage[T2A]{fontenc}
\usepackage[cp1251]{inputenc}
\usepackage[english,russian]{babel}
}%

\makeatletter
\if@twoside%
\else\fi%
\makeatother
\begin{document}

\title{A weak dichotomy below $\Ei\ti\Et$}
\author{Vladimir Kanovei}
\date{\today}

\maketitle

\begin{abstract}
\noindent
We prove that if $\rE$ is an \eqr\  Borel reducible to 
$\Ei\ti\Et$ then either $\rE$ is Borel reducible to the equality 
of countable sets of reals or $\Ei$ is Borel reducible to 
$\rE.$ 
The \lap{either} case admits further strengthening.
\end{abstract}

Let $\dR=\dn.$
Recall that $\Ei$ and $\Et$ are the \eqr s defined on
the set $\rn$ as follows:
\dm
\bay{rcl}
x\Ei y & \text{iff} &
\sus k_0\:\kaz k\ge k_0\:(x(k)=y(k)\,;\\[0.8ex]

x\Et y & \text{iff} &
\kaz k\:(x(k)\Eo y(k))\,;
\eay
\dm
where $\Eo$ is an \eqr\ defined on $\dR$ so that
\dm
\bay{rcl}
a\Eo b & \text{iff} &
\sus n_0\:\kaz n\ge n_0\:(a(n)=b(n)\,.
\eay
\dm
The equivalence $\Et$ is often denoted as $(\Eo){}^\om.$

Kechris and Louveau in \cite{hypsm} and 
Kechris and Hjorth in \cite{hk:nd,hk:rd} proved 
that any Borel \eqr\ $\rE$ satisfying 
$\rE\rebs\Ei,$ resp., $\rE\rebs\Et,$ also 
satisfies the non-strict $\rE\reb\Eo$.
Here $\rebs$ and $\reb$ are resp.\ strict and non-strict
relations of Borel reducibility.
Thus if $\rE$ is an \eqr\ on a Borel set $X$\snos
{We consider only Borel sets in Polish spaces.}
and $\rF$ is an \eqr\ on a Borel set $Y$ then
$\rE\reb\rF$ means that there exists a Borel map
$\vt:X\to Y$ such that 
\dm
{x\rE x'}\leqv{\vt(x)\rF\vt(x')}
\dm
holds for all $x,x'\in X.$
Such a map $\vt$ is called a (Borel) {\it reduction\/} of
$\rE$ to $\rF.$
If both $\rE\reb\rF$ and $\rF\reb\rE$ then they write
$\rE\eqb\rF$ (Borel {\it bi-reducibility\/}),
while $\rE\rebs\rF$ (strict reducibility)
means that $\rE\reb\rF$ but not $\rF\reb\rE.$
See the cited papers \cite{hk:nd,hk:rd} or \eg\
\cite{h,ndir} on various aspects of Borel reducibility in
set theory and mathematics in general.

The abovementioned results give a complete description
of the \dd\reb structure of Borel \eqr s below $\Ei$ and
below $\Et.$
It is then a natural step to investigate the
\dd\reb structure below $\Eit$,
where $\Eit=\Ei\ti\Et$ is the product of $\Ei$ and $\Et,$
that is, an equivalence on $\rnd$ defined so that
for any points $\ang{x,\xi}$ and $\ang{y,\eta}$ in $\rnd,$
$\ang{x,\xi}\Eit\ang{y,\eta}$ if and only if   
$x\Ei y$ and $\xi\Et\eta$.

The intended result would be that the \dd\reb cone below
$\Eit$ includes the cones determined separately
by $\Ei$ and $\Et$, together with the disjoint union of
$\Ei$ and $\Et$
(\ie, the union of $\Ei$ and $\Et$ defined on two disjoint
copies of $\rn$),
$\Eit$ itself, and nothing else.
This is however a long shot.
The following theorem, the main result of this note, can be
considered as a small step in this direction.

\bte
\lam{mt}
Suppose that\/ $\rE$ is a Borel \eqr\ and\/
$\rE\reb{\Eit}$.
Then\/ \bfei\ $\rE$ is Borel reducible
to\/ $\rtd$ \bfor\/ $\Ei\reb\rE$.
\ete 

Recall that the \eqr\ $\rtd,$ known as 
``the equality of countable sets of reals'', is defined 
on $\rn$ so that $x\rtd y$ iff 
$\ens{x(n)}{n\in\dN}=\ens{y(n)}{n\in\dN}.$ 
It is known that $\Et\rebs\rtd$ strictly, and there exist 
many Borel \eqr s $\rE$ satisfying $\rE\rebs\rtd$ but 
incomparable with $\Et:$ 
for instance non-hyperfinite Borel countable ones like 
$\Ey.$
The two cases are incompatible because $\Ei$ is known not 
to be Borel reducible to orbit \eqr s of Polish actions 
(to which class $\rtd$ belongs).

A rather elementary argument reduces Theorem~\ref{mt} to
the following:

\bte
\lam{mt'}
Suppose that\/ $P_0\sq\rnd$ is a Borel set. 
Then\/ \bfei\/ the equivalence\/ 
$\Ejt\res P_0$ is Borel reducible
to\/ $\rtd$ \bfor\/ $\Ei\reb{\Ejt\res P_0}$.
\ete

Indeed suppose that $Z$ (a Borel set) is the domain of
$\rE$,
and $\vt:Z\to\rnd$ is a Borel reduction of $\rE$ to $\Eit$.
Let $f:Z\to \dn=\dR$ be an arbitrary Borel injection.
Define another reduction $\vt':Z\to\rnd$ as follows.
Suppose that $z\in Z$ and $\vt(z)=\ang{x,\xi}\in\rnd.$
Put $\vt'(z)=\ang{x',\xi},$ where $x',$ still a point in
$\rn,$ is related to $x$ so that
$x'(n)=x(n)$ for all $n\ge1$ but $x'(0)=f(z).$
Then obviously $\vt(z)$ and $\vt'(z)$ are
\dd\Eit equivalent for all $z\in Z,$ and hence $\vt'$ is
still a Borel reduction of $\rE$ to $\Eit$.
On the other hand, $\vt'$ is an injection (because so is $f$).
It follows that its full image
$P_0=\ran\vt'=\ens{\vt'(z)}{z\in Z}$ is a Borel set in
$\rnd,$ and $\rE \eqb {\Ejt\res P_0}$.

The remainder of the paper contains the proof of
Theorem~\ref{mt'}. 
The partition in two cases is described in Section~\ref{imp}. 
Naturally assuming that $P_0$ is a lightface $\id11$ set, 
Case~1 is essentially the case when for every element 
$\ang{x,\xi}\in P_0$ 
(note that $x\yi\xi$ are points in $\rn$) 
and every $n$ we have 
$x(n)=F(x\rg n,\xi\rle k,\xi\rg k)$ for some $k,$ 
where $F$ is a $\id11$ function \dd\Et invariant \vrt\ the 
3rd argument. 
It easily follows that then the first projection of the 
equivalence class $\ek{\ang{x,\xi}}{\Eit}\cap P_0$ of 
every point $\ang{x,\xi}\in P_0$ 
is at most countable, leading to the \bfei\ option of 
Theorem~\ref{mt'} in Section~\ref{1:2}.

The results of theorems \ref{mt} and \ref{mt'} 
in their \bfei\ parts can hardly
be viewed as satisfactory because one would expect it in the
form: $\rE$ is Borel reducible to $\Et$.
Thus it is a challenging problem to replace $\rtd$ by 
$\Et$ in the theorems.
Attempts to improve the \bfei\ option, so far rather 
insuccessful, lead us to the following theorem established 
in sections \ref{1:3'} and \ref{1:fin}:

\bte
\lam{mt"}
In the\/ \bfei\ case of Theorem~\ref{mt'} there exist a 
hyperfinite \eqr\ $\rG$ on a Borel set\/ $\po\sq\rnd$ 
such that\/ $\Ejt\res P_0$ is Borel reducible to the 
conjunction of\/ $\rG$ and the \eqr\/ $\Et$ acting on the 2nd 
factor of\/ $\rnd.$\snos
{\label{rs}%
The conjunction as indicated is equal to the least \eqr\ 
$\rF$ on $\po$ which includes $\rG$ and satisfies 
${\xi\Et\eta}\imp{\ang{x,\xi}\rF\ang{y,\eta}}$ 
for all $\ang{x,\xi}$ and $\ang{y,\eta}$ in $\po$.}
\ete

The equivalence $\rG$ as in the theorem will be induced 
by a countable group $\dG$ of homeomorphisms of $\rnd$ 
preserving the second component. 
(That is, if $g\in \dG$ and $g(x,\xi)=\ang{y,\eta}$ then 
$\eta=\xi,$ but $y$ generally speaking depends on both $x$ 
and $\xi$.)
And $\dG$ happens to be even a  
\rit{hyperfinite} group in the sense that it is equal 
to the union of an increasing chain of its finite subgroups. 
Recall that $\Et$ is induced by the product group 
$\dH={\stk{\pwf\dN}\sd}^\dN$ naturally acting in this case 
on the second factor in the product $\rnd.$
And there are further details here that will be presented 
in sections \ref{1:3'} and \ref{1:fin}.

Case 2 is treated in Sections \ref{2:1} through \ref{2:cs}. 
The embedding of $\Ei$ in $\Eit\res P_0$ is obtained by 
approximately the same splitting construction as the one 
introduced in \cite{hypsm} 
(in the version closer to \cite{mlq}).

\punk{Preliminaries: extension of \lap{invariant} functions}
\las{inv}

If $\rE$ is an \eqr\ on a set $X$ 
then, as usual, $\eke x=\ens{y\in X}{y\rE x}$ is the 
\dde\rit{class} of an element $x\in X,$ and 
$\eke Y=\bigcup_{x\in Y}\eke x$   
is the \dde\rit{saturation} of a set $Y\sq X.$
A set $Y\sq X$ is \dde\rit{invariant} if $Y=\eke Y.$ 

The following \lap{invariant} \Sepa\ theorem 
will be used below.

\bpro
[5.1 in \cite{hkl}]
\lam{51}
Assume that\/ $\rE$ is a\/ $\id11$ \eqr\ on a\/ 
$\id11$ set\/ $X\sq\bn.$ 
If\/ $A,C\sq X$ are\/ $\is11$ sets and\/ 
$\eke A\cap\eke C=\pu$ then there exists an\/ 
\dde invariant\/ $\id11$ set\/ $B\sq X$ such that\/ 
$\eke A\sq B$ and\/ $\eke C\cap B=\pu.$
%(Then obviously\/ $\eke A\sq B$ and\/ $\eke C\cap B=\pu.$)
\qed
\epro

Suppose that $f$ is a map defined on a set $Y\sq X.$ 
Say that $f$ is \rit{\dde invariant} 
if $f(x)=f(y)$ for all $x,y\in Y$ satisfying $x\rE y.$ 

\bcor
\lam{51c}
Assume that\/ $\rE$ is a\/ $\id11$ \eqr\ on a\/ 
$\id11$ set\/ $A\sq\bn,$ and\/ $f:B\to\bn$ is an\/ 
\dde invariant\/ $\is11$ function defined on a\/ 
$\is11$ set\/ $B\sq A.$ 
Then there exist an\/ \dde invariant\/ 
%$\id11$ set\/ $Z\sq A$ and an\/ \dde bi-invariant\/ 
$\id11$ function\/ 
$g:A\to\bn$ such that\/ 
%$Y\sq Z$ and\/ 
$f\sq g$. 
%(that is, $g$ extends\/ $f$).
\ecor
\bpf
It obviously suffices to define such a function on 
an \dde invariant $\id11$ set $Z$ such that 
$Y\sq Z\sq A.$
(Indeed then define $g$ to be just a constant on
$A\dif Z$.)
The set 
\dm
P=\ens{\ang{a,x}\in A\ti\bn}
{\kaz b\:\big({(b\in B\land a\rE b)}\imp {x=f(b)}\big)}
\dm
is $\ip11$ and $f\sq P.$ 
Moreover $P$ is \ddf invariant, where 
$\rF$ is defined on $A\ti\bn$ so that 
$\ang{a,x}\rF\ang{a',y}$ iff $a\rE a'$ and $x=y.$ 
Obviously $\ekf f\sq P.$ 
Hence by Proposition~\ref{51} there exists an 
\ddf invariant $\id11$ set $Q$ such that 
$f\sq Q\sq P.$ 
The set
\dm
R=\ens{\ang{a,x}\in Q}
{\kaz y\:({y\ne x}\imp {\ang{a,y}\nin Q}} 
\dm 
is an \ddf invariant $\ip11$ set, 
and in fact a function, satisfying 
$f\sq R.$ 
Applying Proposition~\ref{51} once again we end the proof.
\epf

\punk{An important population of $\is11$ functions}
\las{imp}

Working with elements and subsets of $\rnd$ as the domain
of the \eqr\ $\Eit,$ we'll typically use letters $x,y,z$
to denote
points of the first copy of $\rn$ (where $\Ei$ lives) and
letters $\xi,\eta,\za$ to denote
points of the second copy of $\rn$ (where $\Et$ lives).
Recall that, for $P\sq\rnd,$
\dm
\dom P=\ens{x}{\sus\xi\:(\ang{x,\xi}\in P)}
\quad\text{and}\quad
\ran P=\ens{\xi}{\sus x\:(\ang{x,\xi}\in P)}.
\dm
Points of $\dR=\dn$ will be denoted by $a,b,c$.

Assume that $x\in\rn.$
Let $x\rg n,$ resp., $x\rge n$ denote the restriction of
$x$ (as a map $\dN\to\dR$) to the domain $\iv n\iy,$
resp., $\iry n.$
Thus $x\rg n\in \drb n,$ where ${>}n$ means the
interval ${\iv n\iy},$ and
$x\rge n\in\drbe n,$ where ${\be}n$ means
${\ir n\iy}.$
If $X\sq\rn$ then put
$X\rg n=\ens{x\rg n}{x\in X}$ and
$X\rge n=\ens{x\rge n}{x\in X}$.

The notation connected with $\rl n$ and $\rle n$ is 
understood similarly.

Let $\xi\eqn k\eta$ mean that $\xi\Et\eta$ and
$\xi\rl k=\eta\rl k$
(that is, $\xi(j)=\eta(j)$ for all $j<k$).
This is a Borel equivalence on $\rn.$ 
A set $U\sq\rn$ is \dd{\eqn k}\rit{invariant} if  
$U=\etk U{k},$ where 
$\etk U{k}=\bigcup_{\xi\in U}\ek\xi{\eqn k}$. 

\bdf
\lam{sn}
Let $\xs nk$ denote the set of all 
$\is11$ functions\snos
{A $\is11$ function is a function with a $\is11$ graph.}
$\vpi:U\to\dR\,,$ defined on a $\is11$ set
$U=\dom \vpi\sq \rnn n,$ and 
\dd{\eqn k}\rit{invariant} 
in the sense that if 
$\ang{y,\xi}$ and $\ang{y,\eta}$ belong to
$U$ and $\xi\eqn k\eta$
then $\vpi(y,\xi)=\vpi(y,\eta)$.

%A function $P\in\xs nk$ is 
%\dd{\eqn k}\rit{bi-invariant} iff $D=\dom P$ is a 
%\dd{\eqn k}invariant set. 
Let $\xb nk$ denote the set of all \rit{total}
%\dd{\eqn k}\rit{bi-invariant} 
functions in $\xs nk,$ that is, those defined on 
the whole set $\rnn n$.
\edf 
                           
\ble
\lam{utv1}
If\/ $\vpi\in\xs nk$ then there is a\/ $\id11$ function\/ 
$\psi\in\xb nk$ with\/ $\vpi\sq \psi.$ 
\ele
\bpf
Apply Corollary~\ref{51c}.
\epf

\bdf
\lam{codd}
Let us fix a suitable coding 
system $\sis{\dl e}{e\in E}$ of all
$\id11$ sets $W\sq\dR\ti\rn\ti\dR$ 
(in particular for partial $\id11$ functions 
$\dR\ti\rn\to\dR$), 
where $E\sq\dN$ is a $\ip11$ set, such that there
exist a $\is11$ relation 
$\fSg$ and a $\ip11$ relation $\fPi$ satisfying  
\bus
\label{u1}
\ang{b,\xi,a}\in \dl e\leqv
\fSg(e,b,a,\xi)\leqv \fPi(e,b,a,\xi) 
\eus
whenever $e\in E$ and 
$a,b\in\dR\yt \xi\in\rn$.

Let us fix a $\id11$ sequence of homeomorphisms 
$H_n:\dR\onto\drbe n.$ 
Put
\bus
\left.
\bay{rcl}
\label{u2}
\dk ne&=&
\ens{\ang{H_n(b),\xi,a}}{\ang{b,\xi,a}\in\dl e}
\quad\text{for $e\in E$}\\[\dxii]

T&=&\ens{\ang{e,k}}{e\in E\land \dl e\,
\text{ is a total and \dd{\eqn k}invariant function}}
\eay
\right\}
\eus
Here the totality means that $\dom \dl e=\dR\ti\rn$ 
while the invariance means that 
$\dl e(b,\xi)=\dl e(b,\eta)$ for all $b,\xi,\eta$ 
satisfying $\xi\eqn k\eta.$ 
\edf

Note that if $\ang{e,k}\in T$ then, for any $n,$ 
$\dk ne$ is a function in $\xb nk$, and conversely, every
function in $\xb nk$ has the form $\dk ne$ for a suitable
$e\in E$. 

\bpro
\lam{T} 
$T$ is a\/ $\ip11$ set.
\epro
\bpf
Standard evaluation based on the coding of $\id11$ sets.
\epf

\bcor
\lam{ankp}
%\snos
%{$S^k_n$ is the union of all sets in $\xs nk$.}\/
The sets
\dm
\bay{rcl}
S^k_n &=&\ens{\ang{x,\xi}\in\rnd}
{\sus \vpi \in\xs nk\:(x(n)=\vpi (x\rg n,\xi))}\\[0.8\dxii]
&=&\ens{\ang{x,\xi}\in\rnd}
{\sus \vpi \in\xb nk\:(x(n)=\vpi (x\rg n,\xi))}
\eay
\dm
belong to\/ $\ip11$ uniformly on\/ $n\yi k.$
Therefore the set\/
$\ds=\bigcup_m\bigcap_{n\ge m}\bigcup_k S^k_n$
also belongs to\/ $\ip11$.
\ecor
\bpf
The equality of the two definitions  
follows from Lemma~\ref{utv1}.  
The definability follows from Proposition~\ref{T} 
by standard evaluation.
\epf

%\punk{The cases}
%\las{ca}

Beginning {\ubf the proof of Theorem~\ref{mt'},}
we can \noo\ assume, as usual, that the Borel set 
$P_0$ in the theorem is a lightface $\id11$ set. 

\bde
\item[Case 1:]
$P_0\sq \ds.$
We'll show that
in this case $\Ejt\res P_0$ is 
%\lap{effectively} 
Borel reducible to\/ $\rtd$.

\item[Case 2:]
$P_0\dif \ds\ne\pu.$  
We'll prove that then $\Ei\reb{\Ejt\res P_0}$.
\ede

\punk{Case 1: simplification}
\las{1:1}  

From now on and until the end of Section~\ref{1:2} 
we work under the assumptions of Case~1.
The general strategy is to prove that for any
$\ang{x,\xi}\in P_0$ there exist at most countably many
points $y\in\rn$ such that, for some $\eta,$
$\ang{y,\eta}\in P_0$ and  
$\ang{x,\xi}\Eit\ang{y,\eta},$ and that those points can 
be arranged in countable sequences in a certain controlled 
way.

Our first goal is to somewhat simplify the picture.

\ble
\lam{susl}
There exists a\/ $\id11$ map\/ $\mu:P_0\to\dN$ 
such that for any\/ $\ang{x,\xi}\in P_0$ 
we have\/ 
$\ang{x,\xi}\in\bigcap_{n\ge\mu(x,\xi)}\bigcup_k S^k_n$. 
\ele
\bpf
Apply \Kres\ to the set
\dm
\ens{\ang{\ang{x,\xi},m}\in P_0\ti\dN}
{\kaz n\ge m\:\sus k\:(\ang{x,\xi}\in S^k_n)}\,.
\eqno\qed
\dm
\ePf

Let $\fo=0^\dN\in\dR=\dn$ be the constant $0:$  
$\fo(k)=0\zd\kaz k.$ 
For any $\ang{x,\xi}\in P_0$ put 
$f_\mu(x,\xi)=\fo^{\mu(x,\xi)}\we (x\rge{\mu(x,\xi)}):$ 
that is, we replace by $\fo$ all values $x(n)$ with 
$n<\mu(x,\xi).$ 
Then 
$P'_0=\ens{\ang{f_\mu(x,\xi),\xi}}{\ang{x,\xi}\in P_0}$ 
is a $\is11$ set. 

Put $\ds'=\bigcap_{n}\bigcup_k S^k_n$ 
(a $\ip11$ set by Corollary~\ref{ankp}).

\bcor
\lam{susk}
There is a\/ $\id11$ set\/ $\po$ such that 
$P'_0\sq \po\sq\ds'.$ 
The map 
$\ang{x,\xi}\mto \ang{f_\mu(x,\xi),\xi}$ is a reduction 
of $\Ejt\res P_0$ to $\Ejt\res \po$.
\ecor
\bpf
Obviously $P'_0$ is a subset of the $\ip11$ set 
$\ds'.$ 
It follows that there is a $\id11$ set $\po$ such that 
$P'_0\sq \po\sq\ds'.$ 
To prove the second claim note that ${f_\mu(x,\xi)}\Ei x$ 
for all $\ang{x,\xi}\in P_0$. 
\epf

Let us fix a $\id11$ set $\po$ as indicated. 
By Corollary~\ref{susk} to accomplish Case~1 it suffices 
to get a Borel 
reduction of $\Ejt\res \po$ to\/ $\rtd$.

\ble
\lam{sim}
There exist$:$ a\/ $\id11$ sequence\/ 
$\sis{\kk n}{n\in\dN}$ of natural numbers,
and a\/ $\id11$ system\/ 
$\sis{\xg in}{i,n\in\dN}$ of functions\/
$\xg in\in \xb n{\kk i},$ such that for all\/ 
$\ang{x,\xi}\in \po$ and\/ $n\in\dN$ there is\/ 
$i\in\dN$ satisfying\/ $x(n)=\xg in(x\rg n,\xi)$.
\ele

\bre
\lam{rec}
Recall that by definition every function 
$F \in \xb nk$ 
is invariant in the sense that if 
$\ang{x,\xi}$ and $\ang{x,\eta}$ belong to
$\rnn n,$ $\xi\rl k=\eta\rl k,$ and $\xi\Et\eta,$ 
then $\vpi(x,\xi)=\vpi(x,\eta).$ 
This allows us to sometimes use the notation like
$\xg in(x\rg n,\xi\rl k,\xi\rge k),$ where 
$k=\kk i,$ instead of $\xg in(x\rg n,\xi),$ 
with the understanding that 
$\xg in(x\rg n,\xi\rl k,\xi\rge k)$ is \dd\Et invariant 
in the 3rd argument. 

In these terms,  the final equality of the lemma can be 
re-written as 
$x(n)=\xg in(x\rg n,\xi\rl k,\xi\rge k),$
where $k=\kk i$. 
\ere

\bpf[lemma]
By definition $\po\sq\ds'$ means that for any 
$\ang{x,\xi}\in \po$ and 
$n$ there exists $k$ such that $\ang{x,\xi}\in S^k_n.$
The formula $\ang{x,\xi}\in S^k_n$ takes the form 
\dm
\sus\vpi\in\xb nk\:(x(n)=\vpi(x\rg n,\xi)),
\dm
and further the form
$\sus\ang{e,k}\in T\:(x(n)=\dk ne(x\rg n,\xi)).$ 
It follows that the $\ip11$ set 
\dm
Z=\ens{\ang{\ang{x,\xi,n},\ang{e,k}}\in(P_0\ti\dN)\ti T}
{x(n)=\dk ne(x\rg n,\xi)}
\dm
satisfies $\dom Z=P_0\ti\dN.$ 
Therefore by \Kres\ there is a $\id11$ map 
$\ve:P_0\ti\dN\to T$ such that $x(n)=\dk ne(x\rg n,\xi)$ 
holds for any $\ang{x,\xi}\in P_0$ and $n,$ where 
$\ang{e,k}=\ve(x,\xi,n)$ for some $k.$ 

The range $R=\ran\ve$ of this function is a $\is11$ 
subset of the $\ip11$ set $T.$ 
We conclude that there is a $\id11$ set $B$ such that 
$R\sq B\sq T.$ 
And since $T\sq\dN\ti\dN,$ it follows, by some known 
theorems of effective descriptive set theory, that 
the set $\wE=\dom B=\ens{e}{\sus k\:(\ang{e,k}\in B)}$ 
is $\id11,$ and in addition there exists a $\id11$ map 
$K:\wE\to\dN$ such that $\ang{e,K(e)}\in B$ 
(and $\in T$) for all $e\in\wE$.

And on the other hand it follows from the construction 
that
\bus
\label{u3}
\kaz\ang{x,\xi}\in P_0\:\kaz n\:\sus e\in\wE\:
(x(n)=\dk ne(x\rg n,\xi))\,.
\eus
Let us fix any $\id11$ enumeration 
$\sis{e(i)}{i\in\dN}$ of elements of $\wE.$ 
Put $\xg in=W^{e(i)}_n.$ 
Then the last conclusion of the lemma follows from \eqref{u3}.
Note that the functions $\xg in$ are uniformly $\id11,$ 
$\xg in\in\xb nk$ for some $k,$ in particular, for 
$k=\kk i$, where $\kk i=K(e(i)),$ and $\sis{\kk i}{i\in\dN}$  
is a $\id11$ sequence as well.
\epf

\bbl
\lam{pkf}
Below, we assume that the set $\po$ is chosen as above, 
that is, $\id11$ and $\po\sq\ds',$ while a system of 
functions $\xg in$ and a sequence $\sis{\kk i}{i\in\dN}$ 
of natural numbers 
are chosen accordingly to Lemma~\ref{sim}. 
\ebl

\punk{Case 1: countability of projections of 
equivalence classes}
\las{1:2}

We prove here that in the assumption of Case~1  
the equivalence ${\Ejt}\res \po$ is Borel reducible 
to $\rtd,$ the equality of countable sets of reals. 
%Our next point of interest are sets of the form
The main ingredient of this result will be the 
countability of the sets
\dm
\ikl x\xi=\dom {(\ekit{\ang{x,\xi}}\cap \po)}=
\ens{y\in\rn}{y\Ei x\land \sus \eta\:
(\xi\Et\eta\land\ang{y,\eta}\in \po)},
\dm
where $\ang{x,\xi}\in \po$ --- projections of 
\dd\Eit classes of elements of the set $\po$.

\ble
\lam{ikl}
If\/ $\ang{x,\xi}\in \po$ then\/ $\ikl x\xi\sq\ek x\Ei$ 
and\/ $\ikl x\xi$ is at most countable.
\ele
\bpf 
That $\ikl x\xi\sq\ek x\Ei$ is obvious.
The proof of countability begins with several definitions.
In fact we are going to organize elements of any set of 
the form $\ikl x\xi$ in a countable sequence.

Recall that $\dR=\dn.$
If $u\sq\dN$ and $b\in\dR$ then define $u\ap a\in\dR$ 
so that $(u\ap a)(j)=a(j)$ whenever $j\nin u,$ and 
$(u\ap a)(j)=1-a(j)$ otherwise.

If $f\sq\nn$ and $a\in\dR^k$ 
then define $f\ap a \in\dR^k$ so that 
$(f\ap a)(j)={(\imb fj)}\ap a(j)$ 
for all $j<k,$ where $\imb fj=\ens{m}{\ang{j,m}\in f}.$
Note that $f\ap a$ depends in this case 
only on the restricted set 
$f\res k=\ens{\ang{j,m}\in f}{j<k}.$ 

Put $\Phi=\pwf{\nn}$ and $D=\bigcup_nD_n$, where  
for every $n$:
\dm
D_n=\ens{\ang{a,\vpi}}{a\in\dN^n\land \vpi\in\Phi^n\land
\kaz j<n\:\big(\vpi(j)\sq\kk{a(j)}\ti\dN\big)}.\snom
\dm
(The inclusion $\vpi(j)\sq\kk{a(j)}\ti\dN$ here means 
that the set $\vpi(j)\sq\nn$ satisfies 
$\vpi(j)=\vpi(j)\res{\kk{a(j)}},$ that is, every pair 
$\ang{k,l}\in\vpi(j)$ satisfies $k<\kk{a(j)}$.)

If $\ang{a,\vpi}\in D_n$ and 
$\ang{x,\xi}\in\rnd$ then we define 
$y=\pip x\xi a\vpi\in\rn$ as follows: 
$y=\ang{b_0,b_1,\dots,b_{n-1}}\we {(x\rge n)},$ 
where the reals $b_m\in\dR$ $(m<n)$ are defined by inverse 
induction so that
\bus
\label{bm}
b_m=\xg {a(m)}m\big(\ang{b_{m+1},b_{m+2},\dots,b_{n-1}}\we 
{(x\rge n)}
\,,\,
\vpi(m)\ap{(\xi\rl{\kk{a(m)}})}\,,\,\xi\rge{\kk{a(m)}}\big).
\eus 
(See Remark~\ref{rec} on notation. 
The element
$
\eta=
{\big(\vpi(m)\ap{(\xi\rl{\kk{a(m)}})}\big)}
\we{(\xi\rge{\kk{a(m)}})}
$
belongs to $\rn$ and satisfies $\eta\Et\xi$ because 
$\vpi(m)$ is a finite set.) 

Put $\pip x\xi\La\La=x$    
($\La$ is the empty sequence). 

Note that by definition the element $y=\pip x\xi a\vpi\in\rn$ 
satisfies $y\rge n=x\rge n$ provided $\ang{a,\vpi}\in D_n$, 
thus in any case $x\Ei \pip x\xi a\vpi.$ 
Thus $\piq x\xi,$ the \rit{trace} of $\ang{x,\xi},$ 
is a countable sequence, that is, 
a function defined on $D=\bigcup_nD_n$, a countable set,  
and the set  
$\ran\piq x\xi=\ens{\pip x\xi a\vpi}{\ang{a,\vpi}\in D}$ 
of all terms of this sequence 
is at most countable and satisfies  
$x=\pip x\xi\La\La\in\ran\piq x\xi\sq\ek x{\Ei}$.

\bcl
\lam{icl} 
Suppose that\/ $\ang{x,\xi}\in \po.$
Then\/ $\ikl x\xi\sq \ran\piq x\xi$ --- and hence\/ 
$\ikl x\xi$ is at most countable.
More exactly if\/ $y\in \ikl x\xi$ and\/ $y\rge n=x\rge n$
then there is a pair\/ $\ang{a,\vpi}\in D_n$ such 
that\/ $y=\pip x\xi a\vpi$.
\ecl 

We prove the second, more exact part of the claim.
By definition there is $\eta\in\rn$ such that 
$\ang{y,\eta}\in \po$ and $\xi\Et\eta.$ 
Put $b_m=y(m)\zd\kaz m.$ 
Note that for every $m<n$ there is a number $a(m)$ 
such that 
\dm
\bay{rcll}
b_m&=&
\xg {a(m)}m\big(\ang{b_{m+1},\dots,b_{n-1}}\we 
{(y\rge n)}\,,\,\eta\big) &=\\[\dxii]

&=&\xg {a(m)}m\big(\ang{b_{m+1},\dots,b_{n-1}}\we 
{(y\rge n)}
\,,\,
{\eta\rl{\kk{a(m)}}}\,,\,\eta\rge{\kk{a(m)}}\big)& 
\eay
\dm 
for all $m<n$ (see Blanket Agreement~\ref{pkf}), and hence 
\dm
b_m=\xg {a(m)}m\big(\ang{b_{m+1},\dots,b_{n-1}}\we 
{(x\rge n)}
\,,\,
{\eta\rl{\kk{a(m)}}}\,,\,\xi\rge{\kk{a(m)}}\big)  
\dm 
by the invariance of functions $\xg im$ and because 
$x\rge n=y\rge n.$ 
On the other hand, it follows from the assumption 
$\xi\Et\eta$ that for every $m<n$ there is a finite set 
$\vpi(m)\sq\kk{a(m)}\ti\dN$ such that 
$\eta\rl{\kk{a(m)}}=\vpi(m)\ap{(\xi\rl{\kk{a(m)}})}.$  
Then
\dm
b_m=\xg {a(m)}m\big(\ang{b_{m+1},\dots,b_{n-1}}\we 
{(x\rge n)}
\,,\,
{\vpi(m)\ap{(\xi\rl{\kk{a(m)}})}}\,,\,\xi\rge{\kk{a(m)}}\big)  
\dm 
for every $m<n,$ that is, $y=\pip x\xi a\vpi$, as required.
\epF{Claim and Lemma~\ref{ikl}}

The next result reduces the \eqr\  
$\Eit\res\po$ to the equality of sets of the form 
$\ran\piq x\xi$, that is essentially to the \eqr\ 
$\rtd$ of \lap{equality of countable sets of reals}.

\bcor
\lam{c1.}
Suppose that\/ $\ang{x,\xi}$ and\/ $\ang{y,\eta}$ belong 
to\/ $\po.$
Then\/ ${\ang{x,\xi}}\Eit{\ang{y,\eta}}$
holds if and only if\/ $\xi\Et\eta$ and\/ 
$\ran{\piq x\xi}=\ran{\piq y\eta}$.
\ecor
\bpf
The \lap{if} direction is rather easy. 
If $\xi\Et\eta$ and $\ran\piq y\eta= \ran\piq x\xi$ 
then $x\Ei y$ because $\ran\piq y\eta\sq\ek y\Ei$ and 
$\ran\piq x\xi\sq\ek x\Ei$ by Lemma~\ref{ikl}. 

To prove the converse 
suppose that $\ang{x,\xi}\Eit\ang{y,\eta}.$
Then $\xi\Et\eta$, of course. 
Furthermore, $x\Ei y,$ therefore $x\rge n=y\rge n$ for 
an appropriate $n.$ 
Let us prove that $\ran\piq y\eta= \ran\piq x\xi.$ 
First of all, by definition we have $y\in \ikl x\xi,$ 
and hence (see the proof of Claim~\ref{icl}) 
there exists a pair $\ang{a,\vpi}\in D_n$ such 
that $y=\pip x\xi a\vpi.$ 

Now, let us establish $\ran\piq x\xi=\ran\piq y\xi$ 
(with one and the same $\xi$).
Suppose that $z\in \ran\piq x\xi,$ that is, 
$z= \pip x\xi b\psi$ for a pair $\ang{b,\psi}\in D_m$ 
for some $m.$ 
If $m\ge n$ then obviously 
$z=\pip x\xi b\psi=\pip y\xi b\psi,$ and hence 
(as $x\rge n=y\rge n$) $z\in \ran\piq y\xi.$ 
If $m<n$ then $z=\pip x\xi b\psi=\pip y\xi{a'}{\vpi'},$
where $a'=b\we{( a\rge{m})}$ and 
$\vpi'=\psi\we{(\vpi\rge{m})},$
and once again $z\in \ran\piq y\xi.$ 
Thus $\ran\piq x\xi\sq\ran\piq y\xi.$
The proof of the inverse inclusion 
$\ran\piq y\xi\sq\ran\piq x\xi$ is similar.

Thus $\ran\piq y\xi= \ran\piq x\xi.$ 
It remains to prove $\ran\piq y\eta= \ran\piq y\xi$ 
for all $y\yi\xi\yi\eta$ such that $\xi\Et\eta.$  
Here we need another block of definitions. 

Let $\dH$ be the set of all sets $\da\sq\nn$ 
such that $\imb{\da}j=\ens{m}{\ang{j,m}\in\da}$ 
is finite for all $j\in\dN.$ 
For instance if $\xi,\eta\in\rn$ satisfy $\xi\Et\eta$ 
then the set  
\dm
\xda\xi\eta=\ens{\ang{j,m}}{\xi(j)(m)\ne\eta(j)(m)} 
\dm 
belongs to $\dH.$  
The operation of symmetric difference $\sd$ converts 
$\dH$ into a Polish group equal to the 
product group ${\stk{\pwf\dN}{\sd}}^\dN.$ 

If $n\in\dN\yt\ang{ a,\vpi}\in D_n,$
and $\da\in\dH$ then we define a sequence 
$\vpi'=H_{\da}^ a(\vpi)\in\Phi^n$ so that  
$\vpi'(m)={(\da\res\kk{ a(m)})}\sd\vpi(m)$ for every 
$m<n.$\snos
{Recall that $\da\res k=\ens{\ang{j,i}\in\da}{j<k}$.}
Then the pair $\ang{ a,H_{\da}^ a(\vpi)}$ 
obviously still belongs to $D_n$ and 
$H_{\da}^ a(H_{\da}^ a(\vpi))=\vpi$.

Coming back to a triple of $y\yi\xi\yi\eta\in\rn$ 
such that $\xi\Et\eta,$ let $\da=\xda{\xi}{\eta}.$ 
A routine verification shows that 
$\pip y\eta a\vpi= \pip y\xi a{H_{\da}^ a(\vpi)}$ for 
all $\ang{a,\vpi}\in D.$ 
It follows that $\ran\piq y\eta= \ran\piq y\xi,$ as 
required.  
\epf

\bcor
\lam{c1,}
The restricted relation\/ $\Eit\res\po$ 
is Borel reducible to\/ $\rtd$.
\ecor
\bpf
Since all $\piq x\xi$ are countable sequences of 
reals, the equality $\ran\piq y\eta= \ran\piq x\xi$ of 
Corollary~\ref{c1.} is Borel reducible to $\rtd.$ 
Thus $\Eit\res\po$ is Borel reducible to $\Et\ti\rtd$ 
by Corollary~\ref{c1.}. 
However it is known that $\Et$ is Borel reducible 
to $\rtd$, and so does $\rtd\ti\rtd$. 
\epf

\qeDD{Case 1 of Theorem~\ref{mt'}}

\punk{%
Case 1: a more elementary (?) transformation group}
\las{1:3'}  

Here we begin the proof of Theorem~\ref{mt"}.
Our plan is to define a countable group $\dG$ of 
homeomorphisms of $\rnd$ such that the induced \eqr\ 
$\rG$ satisfies Theorem~\ref{mt"}. 
We continue to argue under the assumptions of Case~1.

First of all let us define the basic domain of transformations, 
\dm
\bdo= \ens{\ang{x,\xi}\in\rnd}
{\kaz n\:\sus\ang{a,\vpi}\in D_n\:(x=\pip x\xi a\vpi)}.
\dm
This is a closed subset of $\rnd.$   
Applying Claim~\ref{icl} with $y=x$ we obtain

\bcor
\lam{popi}
$\po\sq\bdo$.\qed
\ecor

Suppose that pairs $\ang{a,\vpi}$ and $\ang{b,\psi}$ 
belong to $D_n$ for one and the same $n,$ and 
$\ang{x,\xi}\in\rnd.$ 
We define $\gr a\vpi{b}{\psi}x\xi=\ang{y,\xi}\in\rnd$ 
so that
\dm
y=\left\{
\bay{rcl}
\pip x\xi{b}{\psi} &\text{whenever}&
x=\pip x\xi{a}{\vpi}\\[1.2\dxii]
 
\pip x\xi{a}{\vpi} &\text{whenever}&
x=\pip x\xi{b}{\psi}\\[1.2\dxii]
 
x &\text{whenever}&
\pip x\xi{a}{\vpi}\ne x\ne\pip x\xi{b}{\psi}
\eay
\right.
\dm 
Note that if $\pip x\xi{a}{\vpi}=x=\pip x\xi{b}{\psi}$
then still $y=x$ by either of the two first cases of the 
definition. 
And in any case $y\rge n=x\rge n$ provided 
$\ang{a,\vpi}\in D_n$.

%The following is rather clear:

\ble
\lam{sh2*}
Suppose that\/ $n\in\dN$ and pairs\/ 
$\ang{a,\vpi}\yd \ang{b,\psi}$ 
belong to $D_n$. 
Then\/ $\gp a\vpi{b}{\psi}$ is a homeomorphism 
of\/ $\rnd$ onto itself, and\/ 
$\gp a\vpi{b}{\psi}=\gp{b}{\psi} a\vpi$. 

In addition, $\gp a\vpi{b}{\psi}$ is a homeomorphism 
of\/ $\bdo$ onto itself.  
\ele
\bpf
Suppose that $\ang{x,\xi}$ belongs to $\bdo$ and prove 
that so does $\ang{y,\xi}=\gr a\vpi{b}{\psi}x\xi.$
By definition $y$ coincides with one of 
$x\yi\pip x\xi a\vpi\yi \pip x\xi{b}{\psi}.$
So assume that $y=\pip x\xi b\psi.$ 
Consider any $m,$ we have to show that 
$y=\pip y\xi{a'}{\vpi'}$ for some $\ang{a',\vpi'}\in D_m.$
If $m\le n$ then the pair of $a'=b\res m$ and 
$\vpi'=\psi\res m$ obviously works. 
If $m>n$ then take the pair of 
$a'=b\we({b'\rge n})$ 
and $\vpi'=\psi\we({\psi'\rge n})$ 
where $\ang{b',\psi'}\in D_m$ is an arbitrary pair satisfying 
$x=\pip x\xi{b'}{\psi'}$.
\epf

\ble
\lam{rr} 
Suppose that\/ $\ang{x,\xi}\in \bdo.$
Then$:$
\ben
\tenu{{\rm(\roman{enumi})}}
\itla{rr1}
if\/ $\ang{a,\vpi}\yd\ang{b,\psi}\in D_n$  
and\/ 
$\ang{y,\xi}=\gr a\vpi b\psi x\xi$ then\/ 
$\ran\piq x\xi=\ran\piq y\xi\:;$ 

\itla{rr2}
if\/ $y\in\ran\piq x\xi$ then there exist\/ $n$ 
and pairs\/ $\ang{a,\vpi}\yd\ang{b,\psi}\in D_n$ 
such that\/ $\ang{y,\xi}=\gr a\vpi b\psi{x}{\xi}$.
\een
\ele
\bpf
\ref{rr1}
Consider an arbitrary 
$z=\pip x\xi{a'}{\vpi'}\in \ran\piq x\xi,$ where 
$\ang{a',\vpi'}\in D_m$. 
Once again $y$ coincides with one of 
$x\yi\pip x\xi a\vpi\yi \pip x\xi{b}{\psi},$
so assume that $y=\pip x\xi b\psi.$ 
If $m\ge n$ then obviously 
$z=\pip y\xi{a'}{\vpi'}\in\ran\piq y\xi.$ 
If $m<n$ then we have 
$z=\pip y\xi{b'}{\psi'},$ where 
$b'=a'\we{(b\rge m)}$ and 
$\psi'=\vpi'\we{(\psi\rge m)}$.

\ref{rr2}
If $y\in\ran\piq x\xi$ then by 
definition there is a pair $\ang{b,\psi}$ 
in some $D_n$ such that $y=\pip x\xi b\psi.$ 
Then by the way $x\rge n=y\rge n.$ 
As $\ang{x,\xi}\in\bdo,$  
there is a pair $\ang{a,\vpi}\in D_n$ 
such that $x=\pip x\xi a\vpi.$
Then $\ang{y,\xi}=\gr a\vpi b\psi{x}{\xi}$. 
\epf

Let $\dG$ denote the group of all finite superpositions 
of maps of the form $\gp a\vpi{b}{\psi}$, where 
$\ang{a,\vpi}\yd \ang{b,\psi}$ belong to one and the 
same set $D_n$ as in the lemma.
Thus $\dG$ is a countable group of homeomorphisms of 
$\rnd.$ 
(We'll prove that $\dG$ is even an increasing union of 
its finite subgroups!)
Note that a superposition of the form 
$\gp{ a'}{\vpi'}{ a''}{\vpi''}\circ\gp a\vpi{ a'}{\vpi'}$
does not necessarily coincide with 
$\gp{ a''}{\vpi''}{ a}{\vpi}.$

We are going to prove that the equivalence relation $\rG$ 
induced by $\dG$ on $\bdo$ satisfies Theorem~\ref{mt"}. 
To be more exact, $\rG$ is defined on $\bdo$ so that 
$\ang{x,\xi}\rG\ang{y,\eta}$ iff there exists a homeomorphism 
$g\in\dG$ such that $g(x,\xi)=\ang{y,\eta}.$ 
Note that then by definition $\eta=\xi$.

The hyperfiniteness $\rG$ will be established in the next 
Section. 
Now let us study relations between $\dG$ and $\dH,$ the 
other involved group introduced in the proof 
of Corollary~\ref{c1.}. 
For any $\da\in\dH$ define a homeomorphism $H_\da$ of 
$\rnd$ so that $H_\da(x,\xi)=\ang{x,\eta},$ where simply 
$\eta=\da\sd\xi$ in the sense that 
\dm
\eta(m,j)=
\left\{
\bay{rcl}
\xi(m,j) &\text{whenever}& \ang{m,j}\nin\da\\[1\dxii]

1-\xi(m,j) &\text{whenever}& \ang{m,j}\in\da
\eay
\right.
\dm
(Then obviously $\da=\xda\xi\eta$.)
If $\ga,\da\in\dH$ then the superposition 
$H_\da\circ H_\ga$ coincides with $H_{\ga\sd\da}$, where 
$\sd$ is the symmetric difference, as usual.

Transformations of the form $\gp a\vpi{b}{\psi}$ do not 
commute with those of the form $H_\da,$ yet there exists 
a convenient law of commutation:

\ble
\lam{comm}
Suppose that\/ $n\in\dN$ and pairs\/ $\ang{ a,\vpi}$ and\/ 
$\ang{ b,\psi}$ belong to\/ $D_n$, and\/ $\da\in\dH.$ 
Then the superposition\/ $\gp a\vpi b\psi\circ H_\da$ 
coincides with\/  
$H_\da\circ \gp{ a}{\vpi'}{ b}{\psi'}$, 
where\/ $\vpi'=H_\da^ a(\vpi)$ and\/ 
$\psi'=H_\da^ b(\psi)$.
\ele
\bpf
A routine argument is left for the reader.
\epf

Let us consider the group $\dS$ 
of all homeomorphisms $s:\rnd\to\rnd$ of the form
\bus
\label{S}
s= H_\da\circ g_{\ell-1}\circ g_{\ell-2}\dots
\circ g_{1}\circ g_{0}\,,
\eus
where $\ell\in\dN\yt\da\in\dH,$ and each $g_i$ is a 
homeomorphism of $\rnd$ 
of the form $\gp{a_i}{\vpi_i}{b_i}{\psi_i}$, where the pairs
$\ang{a_i,\vpi_i}\yd\ang{b_i,\psi_i}$ belong to one and the 
same set $D_n\yd n=n_i$. 
(It follows that 
$g_{\ell-1}\circ g_{\ell-2}\dots
\circ g_{1}\circ g_{0}\in\dG$.)

Lemma~\ref{comm} implies that $\dS$ is really a group under 
the operation of superposition. 
For instance if 
$g=\gp a\vpi b\psi$ and $g_1$ belong to $\dG$ 
(and $\ang{a,\vpi}\yd\ang{ b,\psi}$ belong to one and 
the same $D_n$) then the superposition
$
H_\da\circ g\circ H_{\da_1}\circ g_1
$
coincides with
$
H_\da\circ H_{\da_1}\circ g'\circ g_1=
H_{\da\sd\da_1}\circ (g'\circ g_1)\,,
$
where $g'=\gp a{\vpi'}b{\psi'}$ and 
$\vpi'=H_{\da_1}^ a(\vpi)\yt \psi'=H_{\da_1}^ b(\psi)$
as in Lemma~\ref{comm}. 
%Finally $H_\da\circ H_{\da_1}=H_{\da\sd\da_1}$.

Thus $\dS$ seems to be a more complicated group than the 
direct cartesian product of $\dG$ and $\dH$, but on the 
other hand more elementary than the free product 
(of all formal superpositions of elements of both groups). 
A natural action of $\dS$ on $\rnd$ is defined as follows: 
if $s$ is as in \eqref{S} then 
$s\app\ang{x,\xi}=
H_\da(g_{\ell-1}(g_{\ell-2}(\dots g_1(g_0(x,\xi))\dots))).$ 
Let $\rS$ denote the induced orbit equivalence relation. 
One can easily check that both the group $\dS$ and the 
action are Polish. 
On the other hand, $\rS$ is obviously the conjunction of 
$\rG$ and the \eqr\/ $\Et$ acting on the 2nd 
factor of $\rnd,$ in the sense of Theorem~\ref{mt"} and 
footnote~\ref{rs} on page~\pageref{rs}.
Thus the next lemma, together with the result of 
Lemma~\ref{Ghypf} on the hyperfiniteness of $\rG$, 
accomplish  the proof of Theorem~\ref{mt"}.

\ble
\lam{act}
Suppose that\/ $\ang{x,\xi}\yd\ang{y,\eta}\in\po.$ 
Then\/ $\ang{x,\xi}\Eit\ang{y,\eta}$ if and only 
if\/ $\ang{x,\xi}\rS\ang{y,\eta}$.
\ele
\bpf
Suppose that $\ang{x,\xi}\Eit\ang{y,\eta}.$ 
Then $y\in\ran \piq x\xi$ by 
Corollary~\ref{c1.}, and further 
$\ang{x,\xi}\rS\ang{y,\xi}$ by Lemma~\ref{rr}\ref{rr2}. 
It remains to note that $\ang{y,\xi}\rS\ang{y,\eta}$ 
by obvious reasons.

Now suppose that $\ang{x,\xi}\rS\ang{y,\eta}.$ 
Then $\xi\Et\eta,$ and hence by Corollary~\ref{c1,} 
it suffices to prove that $\ran\piq x\xi=\ran\piq y\eta.$ 
This follows from two observations saying that 
transformations in $\dH$ and in $\dG$ preserve 
$\ran\piq\ast \ast.$  
First, if $\ang{x,\xi}\in\rnd\yt\da\in\dH,$ and 
$\ang{y,\xi}=H_\da(x,\xi)$ then $\piq x\eta$ obviously is 
a permutation of $\piq y\eta,$ and hence 
$\ran\piq x\xi=\ran\piq x\eta.$
Second, if $\ang{x,\xi}\in\rnd,$ pairs 
$\ang{a,\vpi}\yd\ang{b,\psi}$ belong to one and 
the same set $D_n$, and $\ang{y,\xi}=\gr a\vpi b\psi x\xi,$ 
then $\ran\piq x\xi=\ran\piq y\xi$ by Lemma~\ref{rr}.
\epf 

\qeDD{Theorem~\ref{mt"} modulo Lemma~\ref{Ghypf}}

\punk{\boldmath
Case 1: the \lap{hyperfiniteness} of the countable group 
$\dG$}
\las{1:fin} 

Lemma~\ref{act} reduces further study of Case~1 of 
Theorem~\ref{mt'} to properties of the group $\dS$ and 
its Polish actions. 
This is an open topic, and maybe the next result, the 
\lap{hyperfiniteness} of $\dG,$ one of the two components 
of $\dS,$ can lead to a more comprehensive study.  
One might think that $\dG$ is a rather complicated 
countable group, perhaps close to the free group on two 
(or countably many) generators. 
The reality is different:

\ble
\lam{Ghypf}
$\dG$ is the union of an increasing sequence of finite 
subgroups, therefore the induced \eqr\/ $\rG$ is 
hyperfinite.
\ele
\bpf
Let us show that a finite set of \lap{generators} 
$\gp a\vpi{ a'}{\vpi'}$ produces only finitely many 
superpositions --- this obviously implies the lemma. 
Suppose that $m\in\dN,$ and $\ang{ a_i,\vpi_i}\in D_{n(i)}$ 
for all $i<m.$ 
Put $G_{ij}= \gp{ a_i}{\vpi_i}{ a_j}{\vpi_j}$ provided 
$n(i)=n(j),$ and let $G_{ij}$ be the identity otherwise. 
Thus all $G_{ij}$ are homeomorphisms of $\bdo.$  
We are going to prove that the set of all superpositions 
of the form $f_0\circ f_1\circ\dots\circ f_\ell$, where 
$\ell$ is an arbitrary natural number and each of $f_k$ 
is equal to one of $G_{ij}$ 
($i,j$ depend on $k$) 
contains only finitely many really different functions.

Note that if $i,j<m$ and $n(i)<n(j)$ then the pair 
\dm
\ang{ a_i\we{( a_j\rge{n(i)})}\,,\,
\vpi_i\we{(\vpi_j\rge{n(i)})}}
\dm
belongs to $D_{n(j)}.$ 
We can \noo\ assume that every such a pair occurs in the 
list of pairs $\ang{a_i,\vpi_i}\yd i<m$.

Let us associate a pair 
$q(x,\xi)=\ang{u_{x\xi},w_{x\xi}}$
of finite sets 
\dm
\bay{rcl}
u_{x\xi}&=&\ens{i<m}{\pip x\xi{a_i}{\vpi_i}=x},
\quad\text{and}\\[1\dxii]

w_{x\xi}&=&\ens{\ang{i,j}}
{i,j<m\,\land\, 
%n(i)=n(j)\,\land\, 
\pip x\xi{a_i}{\vpi_i}=\pip x\xi{a_j}{\vpi_j}}
\eay
\dm
with every point $\ang{x,\xi}\in\bdo.$ 
Put $Q=\pws m\ti\pws{m\ti m},$
a (finite) set including all possible values of $q(\pi)$.

\bcl
\lam{cl}
For every\/ $q=\ang{u,w}\in Q$ and\/ $i,j<m$ 
there exists\/ $\tq=\ang{\zu,\tiw}\in Q$ 
such that\/ 
$q(G_{ij}(x,\xi))=\tq$ 
for all\/ $\ang{x,\xi}\in\bdo$ with\/ $q(x,\xi)=q.$
\ecl
\bpf[Claim]
We can assume that $i\ne j$ and $n(i)=n(j)$ since otherwise 
$G_{ij}(x,\xi)=\ang{x,\xi},$ and hence $\tq=q$ works. 
By the same reason we can \noo\ assume that either 
$i\in u\land j\nin u$ or $i\nin u\land j\in u.$
Let say $i\in u\land j\nin u,$ that is, 
%$\hp{a_i}{\vpi_i}(\pi)=\pi\ne\hp{ a_j}{\vpi_j}(\pi)$.
$\pip x\xi{a_i}{\vpi_i}=x\ne\pip x\xi{a_j}{\vpi_j}$.
Then by definition the element 
$\ang{y,\xi}=G_{ij}(x,\xi)=
\gr{a_i}{\vpi_i}{a_j}{\vpi_j}x\xi$
coincides with $\ang{\pip x\xi{ a_j}{\vpi_j},\xi}.$
Let us compute $\tq=q(y,\xi).$ 

Consider an arbitrary $k<m.$ 
To figure out whether $k\in \zu=u_{y\xi}$ we have to 
determine whether $\pip y\xi{a_k}{\vpi_k}=y$ holds.
If $n(k)\ge n(i)=n(j)$ then obviously 
$\pip y\xi{a_k}{\vpi_k}=\pip x\xi{a_k}{\vpi_k},$
and hence $\pip y\xi{a_k}{\vpi_k}=y$ iff 
$\ang{j,k}\in w.$ 
Suppose that $n(k)< n(i)=n(j).$
Then 
\dm
\pip y\xi{a_k}{\vpi_k}=
\pip{\pip y\xi{a_j}{\vpi_j}}\xi{a_k}{\vpi_k}=
\pip y\xi{b}{\psi}\,,
\dm
where the pair 
$\ang{b,\psi}=
\ang{ a_k\we{( a_j\rge{n(k)})}\,,\,
\vpi_k\we{(\vpi_j\rge{n(k)})}}$ 
is equal to one of the pairs 
$\ang{ a_\nu,\vpi_\nu}\yd \nu<m$ 
(and then $n(\nu)=n(i)=n(j)$).
Thus $\pip y\xi{a_k}{\vpi_k}=y$ iff 
$\pip x\xi{a_\nu}{\vpi_\nu}=\pip x\xi{ a_j}{\vpi_j}$ 
iff $\ang{j,\nu}\in w.$ 

Now consider arbitrary numbers $k,k'<m.$ 
To figure out whether $\ang{k,k'}\in \tiw=w_{y\xi}$ we have to 
determine whether 
$\pip y\xi{a_k}{\vpi_k}=\pip y\xi{ a_{k'}}{\vpi_{k'}}$ 
holds.
As above in the first part of the proof of the claim, there 
exist indices $\nu,\nu'<m$ 
(that depend on $q(\pi)=\ang{u,v}$ but not directly on 
$\ang{x,\xi}$) 
such that 
$\pip y\xi{a_k}{\vpi_k}=\pip x\xi{a_{\nu}}{\vpi_{\nu}}$ 
and 
$\pip y\xi{a_{k'}}{\vpi_{k'}}=
\pip x\xi{a_{\nu'}}{\vpi_{\nu'}}.$ 
And then the equality
$\pip y\xi{a_k}{\vpi_k}=\pip y\xi{a_{k'}}{\vpi_{k'}}$ 
is equivalent to $\ang{\nu,\nu'}\in w$.
\epF{Claim} 

Come back to the proof of Lemma~\ref{Ghypf}.

Consider any $q=\ang{u,w}\in Q.$ 
Then $\bdo_q=\ens{\ang{x,\xi}\in\bdo}{q(x,\xi)=q}$ is a Borel 
subset of $\bdo.$ 
%As $Q$ is a finite set, 
It follows from the claim that for every 
superposition of the form 
$f=f_0\circ f_1\circ\dots\circ f_\ell$, where each of $f_k$ 
is equal to one of $G_{ij}$ 
($i,j$ depend on $k$) 
there exists a sequence 
$k_0,k_1,\dots,k_\ell$ of numbers $k_i<m$ such that 
\dm
f(x,\xi)\:=\:
\big(
\hp{a_{k_0}}{\vpi_{k_0}}\circ
\hp{a_{k_1}}{\vpi_{k_1}}\circ\dots\circ
\hp{a_{k_\ell}}{\vpi_{k_\ell}}
\big)(x,\xi)
\dm
for all $\ang{x,\xi}\in\bdo_q$, where $\hp a\vpi$ 
is a map of $\bdo\to\bdo$  defined 
so that $\hp a\vpi(x,\xi)=\ang{\pip x\xi a\vpi,\xi}$ 
for all $\ang{x,\xi}\in\rnd.$  
In other words  
$f=f_0\circ \dots\circ f_\ell$ coincides with 
the superposition 
$\hp{a_{k_0}}{\vpi_{k_0}}\circ\dots\circ
\hp{a_{k_\ell}}{\vpi_{k_\ell}}$ 
on $\bdo_q$.

Note finally that if 
$\ang{ a,\vpi}\in D_n\yt \ang{ b,\psi}\in D_{n'}$, 
and $n'\le n$ then 
%$\hq a\vpi{\hq b\psi\pi}=\hq a\vpi\pi$ 
$\hp a\vpi(\hp b\psi(x,\xi))=\hp a\vpi(x,\xi)$
for all 
$\ang{x,\xi}\in\bdo.$ 
It follows that the superposition 
$\hp{ a_{k_0}}{\vpi_{k_0}}\circ\dots\circ
\hp{ a_{k_\ell}}{\vpi_{k_\ell}}$ 
will not change as a function if we remove all 
factors $\hp{ a_{k_i}}{\vpi_{k_i}}$ such that 
$n(k_i)\le n(k_j)$ for some $j<i.$ 
The remaining superposition obviously contains at most 
$n=\tmax_{i<m}n(i)$ terms, and hence there exist only 
finitely many superpositions of such a reduced form.

As $Q$ itself is finite, this ends the proof of the lemma.
\epF{Lemma~\ref{Ghypf}}

\qeDD{Theorem~\ref{mt"}}

\punk{Case 2}
\las{2:1}

Then the $\is11$ set $R=P_0\cap\dq,$
where $\dq=\dn\dif \ds$ is the chaotic domain,
\index{zzH@$\dq$}
is non-empty.
Our goal will be to prove that 
$\Ei\reb {\Ejt\res R}$ in this case.
The embedding $\vt:\rn\to R$ will have the property that
any two elements $\ang{x,\xi}$ and $\ang{x',\xi'}$ in
the range $\ran\vt\sq R$ satisfy $\xi\Et\xi',$ so that the
\dd{\xi'}component in the range of $\vt$ is trivial.
And as far as the \dd xcomponent is concerned, the embedding
will resemble the embedding defined in Case~1 of the proof
of the 1st dichotomy theorem in \cite{hypsm}
(see also \cite[Ch.~8]{k-var}).

Recall that sets $S^k_n$ were defined in
Corollary~\ref{ankp}, and by definition
%\pagebreak[0] 
\bus
\label{u6}
\left.
\bay{rcl}
{\ang{x,\xi}\in\dq} &\imp&
\kaz m\:\sus n\ge m\:\kaz k\:
(\ang{x,\xi}\nin S^k_n)\\[1\dxii]
&\imp&
\kaz m\:\sus n\ge m\:\kaz k\:\kaz\vpi\in\xs nk\:
\big(x(n)\ne\vpi(x\rg n,\xi)\big) 
\eay
\right\}.
\eus
in Case~2.
Prove a couple of related technical lemmas.

\ble
\lam{nki}
Each set\/ $S^k_n$ is invariant in the following 
sense$:$ 
if\/ $\ang{x,\xi}\in S^k_n$, $\ang{y,\eta}\in \rnd,$
$x\rge n=y\rge n,$ and\/ $\xi\Et \eta$ then\/ 
$\ang{y,\eta}\in S^k_n$.
\ele
\bpf
Otherwise there is a $\id11$ function 
$\vpi \in\xb nk$ such that $y(n)=\vpi(y\rg n,\eta).$ 
Then $x(n)=\vpi(x\rg n,\eta)$ as well because 
$x\rge n=y\rge n.$
We put 
\dm
u_j=\xi(j)\sd \eta(j)=\ens{m}{\xi(j)(m)\ne\eta(j)(m)}
\dm 
for every $j<k,$ these are finite subsets of $\dN.$ 
If $a\in\dn$ and $u\sq\dN$ then define 
$u\app a\in\dn$ so that $(u\app a)(m)=a(m)$ for $m\nin u,$ 
and $(u\app a)(m)=a(m)$ for $m\nin u.$ 
If $\za\in\rn$ then define $f(\za)\in\rn$ so that 
$f(\za)(j)=u_j\app \za(j)$ for $j<k,$ and 
$f(\za)(j)=\za(j)$ for $j\ge k$.

Finally, put $\psi(z,\za)=\vpi(z,f(\za))$ for every 
$\ang{z,\za}\in\rnn n.$ 
The map $\psi$ obviously belongs to $\xb nk$  together 
with $\vpi.$ 
Moreover 
\dm
x(n)=\vpi(x\rg n,\eta)=
\psi(x\rg n,f(\eta))=\psi(x\rg n,\xi)
\dm
because $f(\eta)\rl k=\xi\rl k$, and this contradicts 
to the choice of $\ang{x,\xi}$. 
\epf

The next simple lemma will allow us to split 
$\is11$ sets in $\rnd$.

\ble
\lam{nsk}
If\/ $P\sq\rnd$ is a\/ $\is11$ set and\/ $P\not\sq S^k_n$  
%If\/ $\ang{x,\xi}\nin S^k_n$ then in any\/ $\is11$ set\/ 
%$P\sq\rnd$ containing\/ $\ang{x,\xi}$
then there exist points\/ 
$\ang{x,\xi}$ and\/ $\ang{y,\eta}$ in\/ $P$ with
\dm
y\rg n=x\rg n,\quad \eta\Et\xi,\quad \eta\rl k=\xi\rl k, 
\quad\text{but}\quad y(n)\ne x(n)\,. 
\dm
\ele
\bpf
Otherwise 
$\psi=\ens{\ang{\ang{y\rg n,\eta},y(n)}}{\ang{y,\eta}\in P}$ 
is a map in $\xs nk$, and hence $P\sq S^k_n$, 
%
%However by definition $x(n)=\psi(x\rg n,\xi),$ that is, 
%$\ang{x,\xi}\nin S^k_n,$
contradiction.
\epf

\punk{Case 2: splitting system}
\las{2:s}

We apply a splitting construction, developed in \cite{nwf} 
for the study of \lap{ill}founded Sacks iterations.
Below, $2^n$ will typically denote the set of all
dyadic sequences of length $n,$ and $\bse=\bigcup_n2^n$=
all finite dyadic sequences.

The construction involves a map $\vpi:\dN\to\dN$ assuming 
{\ubf infinitely many} values and each its value 
infinitely many times
(but $\ran\vpi$ may be a proper subset of $\dN$),
another map $\pi:\dN\to\dN,$ 
and, for each $u\in\bse,$ a non-empty 
$\is11$ subset $P_u\sq R=\dq\cap P_0$ ---  
which satisfy a quite long list of properties. 

First of all, if $\vpi$ is already defined at least on 
$\ir0n$ and $u\ne v\in2^n$   
then let 
$
\npi[u,v]=
\tmax\ens{\vpi(\ell)}{\ell<n\land u(\ell)\ne v(\ell)}.
%\ens{i\in\dN}{\kaz k<n
%\skl u(k)\ne v(k)\imp i\cl \vpi(k)\skp}\,.
$
And put $\npi[u,u]=-1$ for any $u$.

Now we present the list of requirements 
\ref{z1} -- \ref{zlast}.
\ben
%\tenu{{\rm(\roman{enumi})}}
\cemu
\itla{z1}
if $\vpi(n)\nin\ens{\vpi(\ell)}{\ell<n}$ then 
$\vpi(n)>\vpi(\ell)$ for each $\ell<n$\enuci;

\itla{z1'}
if $u\in2^n$ then 
$P_u\cap {(\bigcup_kS^k_{\vpi(\ell)})}=\pu$ 
for each $\ell<n$\enuci;

\itla{z2}
every $P_u$ is a non-empty $\is11$ subset of $R\cap\dq$\enuci;

\itla{z6}\msur
$P_{u\we i}\sq P_u$ for all $u\in\bse$ and $i=0,1$\enuci;
\een

Two further conditions are related rather to the 
sets $X_u=\dom P_u$.

\ben
%\tenu{{\rm(\roman{enumi})}}
\cemu
\addtocounter{enumi}4
%\itla{z3}
%if $u\in2^n\yt \ang{x,\xi}\in P_u,$ and $k<n,$ then 
%$\vpi(k)\in \nab x$;
%
\itla{z4}
if $u\yi v\in2^n$ then 
$X_u\rg{\npi[u,v]}=X_v\rg{\npi[u,v]}$\enuci;

\itla{z5}
if $u\yi v\in2^n$ then 
$X_u\rge{\npi[u,v]}\cap X_v\rge{\npi[u,v]}=\pu$.
\vyk{\itla{z7}\msur
$\tmax_{u\in2^n}\dia X_u\to0$ as $n\to\iy$ 
(a reasonable 
Polish metric on $\dnd$ is assumed to be fixed);
}
\een

The content of the next condition is some sort of
genericity in the sense of the Gandy -- Harrington
forcing in the space $\rnd,$ that is, the forcing notion
\dm
\dP= \text{ all non-empty $\is11$ subsets of $\rnd$}.
\dm
Let us fix a countable transitive
model $\mm$ of a sufficiently large fragment of $\ZFC.$\snos
{For instance remove the Power Set axiom but add the axiom
saying that for any set $X$ there exists the set of all
countable subsets of $X.$}
For technical reasons, we assume that $\mm$ is an elementary
submodel of the universe \vrt\ all analytic formulas.
Then simple relations between sets in $\dP$ in the universe,
like $P=Q$ or $P\sq Q,$ are adequately reflected as the same
relations between their intersections $P\cap\mm\yt Q\cap\mm$
with the model $\mm.$
In this sense $\dP$ is a forcing notion in
$\mm$.

A set $D\sq\dP$ is {\it open dense\/} iff,
first, for any $P\in\dP$ there is $Q\in D\yt Q\sq P,$
and given sets $P\sq Q\in\dR,$ if $Q$ belongs to $D$ then 
so does $P.$
A set $D\sq\dP$ is {\it coded in $\mm$,\/} iff
the set $\ens{P\cap\mm}{P\in D}$ belongs to $\mm.$ 
There exists at most countably many such sets because $\mm$
is countable. 
Let us fix an enumeration ({\ubf not} in $\mm$)
$\ens{D_n}{n\in\dN}$
of all open dense sets $D\sq\dP$ coded in $\mm$.

The next condition essentially asserts the
\dd\dP genericity of each branch in the splitting
construction over $\mm$.

\ben
%\tenu{{\rm(\roman{enumi})}}
\cemu
\addtocounter{enumi}6
\itla{z8}
for every $n,$ if $u\in 2^{n+1}$ then $P_u\in D_n$.
\een

\bre
\lam{vii}
It follows from \ref{z8} that for any $a\in\dn$ the sequence 
$\sis{P_{a\res n}}{n\in\dN}$ is generic enough for 
the intersection $\bigcap_nP_{a\res n}\ne\pu$ to 
consist of a single point, say $\ang{g(a),\ga(a)},$ 
and for the maps $g,\ga:\dn\to\rnd$ to be continuous.

Note that $g$ is $1-1.$
Indeed if $a\ne b$ belong to $\dn$ then $a(n)\ne b(n)$ for
some $n,$ and hence $\npi[a\res m\,,\,b\res m]\ge\vpi(n)$
for all $m\ge n.$
It follows by \ref{z5} that
$X_{a\res m}\cap X_{b\res m}=\pu$
for $m>n,$ therefore $g(a)\ne g(b)$.
\ere

Our final requirement involves the \dd\xi parts of 
sets $P_u$.
We'll need the following definition. 
Suppose that $\ang{x,\xi}$ and $\ang{y,\eta}$ belong 
to $\rnd,$ $p\in\dN,$ and $s\in\nse\yt\lh s=m$
(the length of $s$).
Define $\ang{x,\xi}\enp ps\ang{y,\eta}$ iff
\dm
\xi\Et\eta\,,\quad 
x\rg{p}=y\rg{p}\,,\quad\text{and}\quad
\xi(k)\sd \eta(k)\sq s(k)\;\;\text{for all}\;\;k<m=\lh s\,,
\dm
where 
$\al\sd\ba=\ens{j}{\al(j)\ne\ba(j)}$ for $\al,\ba\in\dn.$
If $P,Q\sq\rnd$ are arbitrary sets 
then under the same circumstances $P\enp ps Q$
will mean that\pagebreak[0] 
\dm
\kaz \ang{x,\xi}\in P\;\sus\ang{y,\eta}\in Q\;
(\ang{x,\xi}\enp ps\ang{y,\eta})
\quad\text{and \sl vice versa}\,.
\dm
Obviously $\enp ps$ is an equivalence relation.

The following is the last condition:

\ben
%\tenu{{\rm(\roman{enumi})}}
\cemu
\addtocounter{enumi}7
\itla{z'}
\label{zlast}
there exists a map $\pi:\dN\to\dN,$ 
such that 
%we have 
$P_u \enp {\npi[u,v]}{\pi\res n} P_v$
holds for every $n$ and all $u\yi v\in2^n$ \hfill
(and then $X_u\rg{\npi[u,v]}=X_v\rg{\npi[u,v]}$ as in \ref{z4}).
\hfill\,
\een

\punk{Case 2: splitting system implies the reducibility}
\las{2:ssr}

Here we prove that any system of sets $P_u$ and
$X_u=\dom{P_u}$ and 
maps $\vpi\yi\pi$ satisfying \ref{z1} -- \ref{zlast} 
implies Borel reducibility of $\Ei$ to 
${\Ejt}\res R.$ 
This completes Case~2.
The construction of such a splitting system will follow
in the remainder.

Let the maps $g$ and $\ga$ be defined as in
Remark~\ref{vii}.
Put
\dm
W=\ens{\ang{g(a),\ga(a)}}{a\in\dn}.
\dm
%Let us show that $\Ei$ admits an embedding 
%into the restriction of $\Eit$ to $W.$

\ble
\lam{w11}
$W$ is a closed set in\/ $\rnd$ and a function. 
%{\rm(That is, if\/ $\ang{x,\xi}\in W$ and\/ 
%$\ang{x,\eta}\in W$ then\/ $\xi=\eta$.)}
Moreover if\/ $\ang{x,\xi}$ and\/ 
$\ang{y,\eta}$ belong to\/ $W$ then\/ $\xi\Et\eta$.
\ele
\bpf
$W$ is closed as a continuous image of $\dn.$ 
That $W$ is a function follows from the bijectivity 
of $g,$ see Remark~\ref{vii}. 
Finally  
any two $\xi\yi\eta$ as indikated satisfy 
$\xi(k)\sd \eta(k)\sq\pi(k)$ for all $k$ by \ref{z'}.
\epf

Put $X=\dom W.$ 
Thus $W$ is a continuous map $X\to\rn$ by the lemma.

\bcor
\lam{w12}
There exists a Borel reduction of\/ $\Ei\res X$ to\/
$\Eit\res W.$
\ecor
\bpf
As $W$ is a function, we can use the notation $W(x)$ for
$x\in X=\dom W.$
Put $f(x)=\ang{x,W(x)}.$
This is a Borel, even a continuous map $X\to W.$
It remains to establish the equivalence
\bus
\label{ee}
{x\Ei y}\leqv {f(x)\Eit f(y)}\qquad
\text{for all}\quad x,y\in X.
\eus
If $x\Ei y$ then $W(x)\Et W(y)$ by Lemma~\ref{w11}, and hence
easily $f(x)\Eit f(y).$
If $x\Ei y$ fails then obviously $f(x)\Eit f(y)$ fails, too.
\epf

Thus to complete Case~2 it now suffices to define a 
Borel reduction of $\Ei$ to $\Ei\res X.$ 
To get such a reduction consider the set 
$\Phi=\ran\vpi,$ and let $\Phi=\ens{p_m}{m\in\dN}$
in the increasing order; 
that the set 
$\Phi\sq\dN$ is infinite follows from \ref{z1}. 

Suppose that $n\in\dN.$ 
Then $\vpi(n)=p_m$ for some (unique) $m:$ we put 
$\psi(n)=m.$ 
Thus $\psi:\dN\onto\dN$ and the preimage 
$\psi\obr(m)=\vpi\obr(p_m)$ is an infinite subset of 
$\dN$ for any $m.$ 
Define a parallel system of 
sets $Y_u\sq\rn\yt u\in\bse,$ as follows. 
Put $Y_\La=\rn.$ 
Suppose that $Y_u$ has been defined, $u\in2^n.$ 
Put $p=\vpi(n)=p_{\psi(n)}.$ 
Let $K$ be the number of all indices $\ell<n$ still 
satisfying $\vpi(\ell)=p,$ perhaps $K=0.$ 
Put $Y_{u\we i}=\ens{x\in Y_u}{x(p)(K)=i}$ for 
$i=0,1$. 

Each of $Y_u$ is clearly a basic clopen set in $\rn,$ 
and one easily verifies that conditions 
\ref{z6}, \ref{z4}, \ref{z5}   
are satisfied for the sets $Y_u$ and the map $\psi$ 
(instead of $\vpi$ in \ref{z4}, \ref{z5}), in 
particular
\ben
%\tenu{{\mtho$\rm(\roman{enumi}')$}}
\aemu
\addtocounter{enumi}5
\itla{z4'}
if $u\yi v\in2^n$ then 
$Y_u\rg{\nsi[u,v]}=Y_v\rg{\nsi[u,v]}$;

\itla{z5'}
if $u\yi v\in2^n$ then 
$Y_u\rge{\nsi[u,v]}\cap Y_v\rge{\nsi[u,v]}=\pu$;
\een
where 
$\nsi[u,v]=
\tmax\ens{\psi(\ell)}{\ell<n\land u(\ell)\ne v(\ell)}$
(compare with $\npi$ above).

It is clear that for any $a\in\dn$ the intersection 
$\bigcap_nY_{a\res n}=\ans{f(a)}$ 
is a singleton, and the map $f$ is continuous and $1-1.$ 
(We can, of course, define $f$ explicitly: 
$f(a)(p)(K)=a(n),$ where $n\in\dN$ is chosen so that 
$\psi(n)=p$ and there is exactly $K$ numbers $\ell<n$ 
with $\psi(\ell)=p$.) 
Note finally that $\ens{f(a)}{a\in\dn}=\rn$ since 
by definition $Y_{u\we 1}\cup Y_{u\we 0}=Y_u$ for 
all $u$.

We conclude that the map $\vt(x)=g(f\obr(x))$ is a 
continuous map 
(in fact a homeomorphism in this case by compactness) 
$\rn\onto X=\dom W.$ 

\ble
\lam{tel3}
The map\/ $\vt$ is a reduction of\/ $\Ei$ 
to\/ $\Ei\res X,$
and hence\/ $\vt$ witnesses\/ $\Ei\reb{\Ei\res X}$ 
and\/ $\Ei\reb{{\Ejt}\res W}$ by Corollary~\ref{w12}.
\ele
\bpf
It suffices to check that the map $\vt$ satisfies 
the following requirement: 
for each $y\yi y'\in \rn$ and $m$,
\bus
\label{9-0}
y\rge m= y'\rge m\quad\text{iff}\quad
\vt(y)\rge{p_m}=\vt(y')\rge{p_m}
\,.
\eus
To prove \eqref{9-0}
suppose that $y=f(a)$ and $x=g(a)=\vt(y),$ and 
similarly $y'=f(a')$ and $x'=g(a')=\vt(y'),$ where 
$a\yi a'\in\dn.$ 
Suppose that $y\rge m= y'\rge m.$ 
We then have $m>\nsi[a\res n,a'\res n]$ for any $n$ 
by \ref{z5'}. 
It follows, by the definition of $\psi,$ that 
$p_m>\npi[a\res n,a'\res n]$ for any $n,$ hence, 
$X_{a\res n}\rge{p_m}=X_{a'\res n}\rge{p_m}$ 
for any $n$ by \ref{z4}. 
Therefore $x\rge{p_m}=x'\rge{p_m}$ by \ref{z8}, that is, 
the right-hand side of \eqref{9-0}. 
The inverse implication in \eqref{9-0} is proved similarly. 
\epF{Lemma}

It follows that we can now focus on the construction of a
system satisfying \ref{z1} -- \ref{zlast}.
The construction follows in Section~\ref{2:cs}, after
several preliminary lemmas in Sections \ref{2:sh} and
\ref{2:spl}.

\punk{Case 2: how to shrink a splitting system}
\las{2:sh}

Let us prove some results related to 
preservation of condition \ref{z'} under certain 
transformations of shrinking type.
They will be applied in
the construction of a splitting system satisfying
conditions \ref{z1} -- \ref{zlast} of Section~\ref{2:s}.

\ble
\lam{suz}
Suppose that\/ $n\in\dN,$ 
$s\in\nse,$ and a system of\/ $\is11$ sets\/ 
$\pu\ne P_u\sq\rnd\yt u\in2^n,$ satisfies\/ 
$P_u \enp{\npi[u,v]}s P_v$ for all\/ $u,v\in2^n.$  
Assume also that\/ 
$w_0\in2^n,$ and\/ $\pu\ne Q\sq P_{w_0}$ is a\/ 
$\is11$ set.       
Then the system of\/ $\is11$ sets\/
\dm 
P'_u=\ens{\ang{x,\xi}\in P_u} 
{\sus\ang{z,\za}\in Q\,
(\ang{x,\xi}\enp{\npi[u,w_0]}{s}\ang{z,\za})}\,,
\quad u\in2^n,
\dm
still satisfies\/ 
$P'_u \enp{\npi[u,v]}{s} P'_v$ for all\/ $u,v\in2^n,$ 
and $P'_{w_0}=Q.$  
\ele
\bpf
$P'_{w_0}=Q$ holds because $\npi[w_0,w_0]=-1.$ 
%The sets $P'_u$ are as required, via a routine 
%verification.
Let us verify \ref{z'}. 
Suppose that $u,v\in2^n.$
Each one of the three numbers 
$\npi[u,w]\yd\npi[v,w]\yd\npi[u,v]$ is obviously 
not bigger than the largest of the two other numbers. 
This observation leads us to the following three cases.\vom

{\ubf Case a\,:} 
$\npi[u,w_0]=\npi[u,v]\ge\npi[v,w_0].$
%(or the symmetric case).
Consider any $\ang{x,\xi}\in P'_u.$ 
Then by definition there exists $\ang{z,\za}\in Q$ 
with
$\ang{x,\xi}\enp {\npi[u,w_0]}{s}\ang{z,\za}.$ 
Then, as $P_{w_0} \enp {\npi[v,w_0]}{s} P_v$
is assumed by the lemma,  
there is $\ang{y,\eta}\in P_v$ such that 
$\ang{y,\eta}\enp {\npi[v,w_0]}{s}\ang{z,\za}.$ 
Note that $\ang{z,\za}$ witnesses $\ang{y,\eta}\in P'_v$. 
On the other hand, 
$\ang{x,\xi}\enp {\npi[u,v]}{s}\ang{y,\eta}$
because $\npi[u,w_0]=\npi[u,v]\ge\npi[v,w_0].$
Conversely, suppose that $\ang{y,\eta}\in P'_v$.
Then there is $\ang{z,\za}\in Q$ such that 
$\ang{y,\eta}\enp {\npi[v,w_0]}{s}\ang{z,\za}.$ 
Yet $P_{w_0} \enp {\npi[u,w_0]}{s} P_u$, and hence there
exists $\ang{x,\xi}\in P'_u$ with
$\ang{x,\xi}\enp {\npi[u,w_0]}{s}\ang{z,\za}.$ 
Once again we conclude that
$\ang{x,\xi}\enp {\npi[u,v]}{s}\ang{y,\eta}$.\vom

{\ubf Case b\,:}  
$\npi[v,w]=\npi[u,v]\ge\npi[u,w].$
Absolutely similar to Case~a.\vom

{\ubf Case c\,:} 
$\npi[u,w_0]=\npi[v,w_0]\ge\npi[u,v].$
This is a symmetric case, thus it is enough to carry
out only the direction $P'_u\to P'_v$. 
Consider any $\ang{x,\xi}\in P'_u.$ 
As above there is $\ang{z,\za}\in Q$ 
such that 
$\ang{x,\xi}\enp {\npi[u,w_0]}{s} \ang{z,\za}.$ 
On the other hand, as 
$P_u \enp {\npi[u,v]}{s} P_v$,  
there exists a point $\ang{y,\eta}\in P_v$ 
such that 
$\ang{y,\eta}\enp {\npi[u,v]}{s}\ang{x,\xi}.$ 
Note that $\ang{z,\za}$ witnesses 
$\ang{y,\eta}\in P'_v:$ 
indeed by definition we have 
$\ang{y,\eta}\enp{\npi[v,w_0]}{s}\ang{z,\za}.$%
\epf

\bcor
\lam{suz2}
Assume that\/ $n\in\dN,$ 
$s\in\nse,$ and a system of\/ $\is11$ sets\/ 
$\pu\ne P_u\sq\rnd\yt u\in2^n,$ satisfies\/ 
$P_u \enp {\npi[u,v]}s P_v$ for all\/ $u,v\in2^n.$  
Assume also that\/ $\pu\ne W\sq2^n,$ and a\/ $\is11$ set\/
$\pu\ne Q_{w}\sq P_{w}$ is defined for every\/ $w\in W$
so that still\/
$Q_w \enp {\npi[w,w']}s Q_{w'}$ for all\/ $w,w'\in W.$ 
Then the system of\/ $\is11$ sets\/
\dm
P'_u =\ens{\ang{x,\xi}\in P_u}
{\kaz w\in W\:\sus\ang{y,\eta}\in Q_{w}\:
(\ang{x,\xi}\enp{\npi[u,w]}s\ang{y,\eta})}
\dm
still satisfies\/ $P'_u \enp {\npi[u,v]}s P'_v$ for all\/
$u,v\in2^n,$ and\/ $P'_{w}=Q_w$ for all $w\in W.$   
\ecor
\bpf
Apply the transformation of Lemma~\ref{suz} consecutively
for all $w_0\in W$ and the corresponding sets $Q_{w_0}$. 
Note that these transformations do not change the sets 
$Q_{w}$ with $w\in W$ because
$Q_w \enp {\npi[w,w']}s Q_{w'}$ for all\/ $w,w'\in W.$
\epf

\bre
\lam{suzr}
The sets $P'_u$ in Corollary~\ref{suz2} can as well be
defined by
\dm
P'_u =\ens{\ang{x,\xi}\in P_u}
{\sus\ang{y,\eta}\in Q_{w_u}\:
(\ang{x,\xi}\enp{\npi[u,w_u]}s\ang{y,\eta})}
\dm
where, for each $u\in2^n,$ $w_u$ is an element of $W$ such that
the number $\npi[u,w_u]$ is the least of all numbers of the
form $\npi[u,w]\yt w\in W.$
(If there exist several $w\in W$ with the minimal
$\npi[u,w]$ then take the least of them.)
\ere

\punk{Case 2: how to split a splitting system}
\las{2:spl}

Here we consider a different question related to the
construction of systems satisfying
conditions \ref{z1} -- \ref{zlast} of Section~\ref{2:s}. 
Given a system of $\is11$ sets satisfying a \ref{z'}-like
condition, how to shrink the sets so that \ref{z'} is 
preserved and in addition \ref{z5} holds. 
Let us begin with a basic technical question:
given a pair of
$\is11$ sets $P,Q$ satisfying $P\enp ps Q$ for some $p,s,$
how to define a
pair of smaller $\is11$ sets $P'\sq P\yt Q'\sq Q,$ still
satisfying the same condition, but as disjoint as it is
compatible with this condition.

Recall that $\dom P=\ens{x}{\sus\xi\:(\ang{x,\xi}\in P}$
for $P\sq\rnd.$

\ble
\lam{ld}
If\/ $P,Q\sq\rnd$ are non-empty\/ $\is11$ sets,
$p\in\dN\yt s\in\nse,$ $P\enp ps Q,$ and\/
${(P\cup Q)}\cap S^k_p=\pu,$ where\/ $k=\lh s,$
then there exist non-empty \/ $\is11$ sets\/ 
$P'\sq P\yt Q'\sq Q$ such that still\/
$P'\enp ps Q'$ but in addition\/
$(\dom P')\rge p\cap(\dom Q')\rge p=\pu$.
\ele

Note that $P\enp sp Q$ implies $(\dom P)\rg p=(\dom Q)\rg p$.

\bpf
It follows from Lemma~\ref{nsk} that there exist points  
$\ang{x_0,\xi_0}$ and\/ $\ang{x_1,\xi_1}$ in $P$ such that
$\ang{x_0,\xi_0}\enp ps \ang{x_1,\xi_1}$
but $x_1(p)\ne x_0(p).$
Then there exists a number $j$ such that, say,
$x_1(p)(j)=1\ne0= x_0(p)(j).$
On the other hand, there exists $\ang{y_0,\eta_0}\in Q$
such that $\ang{x_i,\xi_i}\enp ps \ang{y_0,\eta_0}$ for
$i=0,1.$
Then $y_0(p)(j)\ne x_i(p)(j)$ for one of $i=0,1.$ 
Let say $y_0(p)(j)=0\ne 1=x_0(p)(j).$
Then the $\is11$ sets
\dm
\bay{rcl}
P'
&\!\!\!=\!\!\!&
\ens{\ang{x,\xi}\in P}
{\sus\ang{y,\eta}\in Q\:\big(
%\!\!\!&\!\!\!\big(\!\!\!&\!\!\!
x(p)(j)=1 \land y(p)(j)=0\land
\ang{x,\xi}\enp {p}{s}\ang{y,\eta}\big)};\\[1.5\dxii]
Q'
&\!\!\!=\!\!\!&
\ens{\ang{y,\eta}\in Q}
{\sus\ang{x,\xi}\in P \:\big(
x(p)(j)=1 \land y(p)(j)=0\land
\ang{x,\xi}\enp {p}{s}\ang{y,\eta}\big)}
\eay
\dm
are $\is11$ and non-empty 
(contain resp.\ $\ang{x_0,\xi_0}$ and $\ang{y_0,\eta_0}$), 
and they satisfy $P'\enp{p}{s}Q'$, but
$(\dom P')\rge p\cap(\dom Q')\rge p=\pu$ 
because $y(p)(j)=0\ne 1=x(p)(j)$ whenever
$\ang{x,\xi}\in P'$ and $\ang{y,\eta}\in Q'.$
\epf

\bcor
\lam{ldc}
Assume that\/ $n\in\dN\yt s\in\nse,$ and a system of\/
$\is11$ sets\/ $\pu\ne P_u\sq\rnd\yt u\in2^n,$ satisfies\/ 
$P_u \enp {\npi[u,v]}s P_v$ for all\/ $u,v\in2^n.$  
Then there exists a system of\/ $\is11$ sets\/
$\pu\ne P'_u\sq P_u\yt u\in2^n,$ such that still\/
$P'_u \enp {\npi[u,v]}s P_v$, and
in addition\/
$(\dom P'_u)\rge{\npi[u,v]}\cap(\dom P'_v)\rge{\npi[u,v]}=\pu,$
for all\/ $u\ne v\in2^n.$ 
\ecor
\bpf
Consider any pair of $u_0\ne v_0$ in $2^n.$
Apply Lemma~\ref{ld} for the sets $P=P_{u_0}$ and $Q=P_{v_0}$
and $p= \npi[u_0,v_0].$
Let $P'$ and $Q'$ be the $\is11$ sets obtained,
in particular $P'\enp {\npi[u_0,v_0]}s Q'$ and
$(\dom P')\rge{\npi[u_0,v_0]}\cap
(\dom Q')\rge{\npi[u_0,v_0]}=\pu$.
Then by Corollary~\ref{suz2} there is a system of $\is11$
sets $\pu\ne P'_u\sq P_u$ such that still  
$P'_u \enp {\npi[u,v]}s P'_v$ for all
$u,v\in2^n,$ and $P_{u_0}=P'\yt P_{v_0}=Q'$ --- and hence
\dm
(\dom P'_{u_0})\rge{\npi[u_0,v_0]}\cap
(\dom P'_{v_0})\rge{\npi[u_0,v_0]}=\pu.
\dm
Take any other pair of $u_1\ne v_1$ in $2^n$ and transform
the system of sets $P'_u$ the same way.
Iterate this construction sufficient (finite) number of steps.
\epf

\punk{Case 2: the construction of a splitting system}
\las{2:cs}

We continue the proof of Theorem~\ref{mt'} -- Case~2.
Recall that $R=P_0\cap\dq$ is a $\is11$ set. 
By Lemma~\ref{tel3}, it suffices to define functions
$\vpi$ and $\pi$ and a system of 
$\is11$ sets $P_u\sq R$ together satisfying conditions 
\ref{z1} -- \ref{zlast}.
The construction of such a system will go on by induction
on $n.$
That is, at any step $n$ the sets $P_u$ with $u\in2^n,$ as
well as the values of $\vpi(k)$ and $\pi(k)$ with $k<n,$
will be defined.

%As \lap{$\glu x=\piy$} is a $\is11$ relation,
%$R'= \ens{x\in D}{\glu x=\piy}$ is still a $\is11$ set.
%Take any $x_0\in D$ and put $\vpi(0)$ to be any 
%$j\in\nab{x_0}.$ 

For $n=0,$ we put $P_\La=R.$
($\La\in2^0$ is the only sequence of length $0$.)

Suppose that sets $P_u\sq R$ with $u\in2^n,$ and also 
all values $\vpi(\ell)\yt\ell<n,$ and $\pi(k)\yt k<n,$ 
have been defined and satisfy the applicable part of 
\ref{z1} -- \ref{zlast}.
The content of the inductive step $n\mto n+1$ will consist
in definition of $\vpi(n)\yt \pi(n),$ and sets
$P_{u\we i}$ with $u\we i\in 2^{n+1},$ that is, $u\in2^n$
(a dyadic sequence of length $n$)
and $i=0,1.$ 
This goes on in four steps A\yi B\yi C\yi D.
%\vom 

\ppu
{\boldmath Step A: definition of $\vpi(n)$}

Suppose that, in the order of increase, 
\dm
\ens{\vpi(\ell)}{\ell<n}=
\ans{p_0<\dots<p_m}\,.
\dm
%in the increasing order. 
For $j\le m,$ let $K_j$ be the number of all 
$\ell<n$ with $\vpi(\ell)= p_j$.\vom 

{\sl Case A\/}:  
$K_j\ge m$ 
%(then actually $N_p=m$) 
for all $j\le m.$ 
Then consider any $u_0\in2^n$ and an arbitrary point 
$\ang{x_0,\xi_0}\in P_{u_0}$. 
Note that by \eqref{u6} of Section~\ref{2:1} 
there is a number $p>\tmax_{\ell<n}{\vpi(\ell)}$ 
such that $\ang{x_0,\xi_0}\nin \bigcup_kS^k_p$.
Put $\vpi(n)=p$. 
%We have something else to do in this case. 
%Put $X_u=\dom P_u$, as above.

We claim that the sets 
$P'_u=P_u\dif\bigcup_kS^k_{\vpi(n)}$ 
still satisfy condition \ref{z'} 
(and then \ref{z4} for $X'_u=\dom P'_u$).
Indeed suppose that $u,v\in2^n$ and 
$\ang{x,\xi}\in P'_u.$ 
Then $\ang{x,\xi}\in P_u,$ and hence there is a 
point $\ang{y,\eta}\in P_v$ such that 
$\ang{x,\xi}\enq {\npi[u,v]}\pi n\ang{y,\eta}.$
%
%$\xi\Et\eta,$ 
%$x\rg{\npi[u,v]}=y\rg{\npi[u,v]},$ and 
%$\xi(k)\sd \eta(k)\sq \pi(k)$ for all $k<n.$  
It remains to show that 
$\ang{y,\eta}\nin \bigcup_kS^k_{\vpi(n)}$. 
Suppose towards the contrary that 
$\ang{y,\eta}\in S^k_{\vpi(n)}$ for some $k.$ 
By definition ${\vpi(n)}>\npi[u,v],$ therefore  
$x\rge{{\vpi(n)}}=y\rge{{\vpi(n)}}$. 
It follows that $\ang{x,\xi}\in S^k_{\vpi(n)}$ 
by Lemma~\ref{nki}, 
contradiction.\vom

{\sl Case B\/}.  
If some numbers $K_j$ are $<m$ then  
choose $\vpi(n)$ among $p_j$ with the least $K_j$, and 
among them take the least one. 
Thus $\vpi(n)=\vpi(\ell)$ for some $\ell<n.$ 
It follows that in this case 
$P_u\cap{(\bigcup_kS^k_{\vpi(n)})}=\pu$ 
for all $u\in2^n$ by the inductive assumption of 
\ref{z1'}. 
Put $P'_u=P_u$.\vom

Note that this manner of choice of $\vpi(n)$ implies 
\ref{z1}, \ref{z1'} and also implies that $\vpi$ takes 
infinitely many values and takes each its value 
infinitely many times. 
In addition, the construction given above proves:

\ble
\lam{st1}
There exists a system of\/ $\is11$ sets\/
$\pu\ne P'_u\sq P_u$ satisfying\/ \ref{z'} and\/ 
$P'_u\cap{(\bigcup_kS^k_{\vpi(n)})}=\pu$ 
for all\/ $u\in2^n.$\qed
\ele

%Let us fix such a system of sets $P'_u\sq P_u$.\vom

\ppu
{\boldmath Step B: definition of $\pi(n)$}

We work with the sets $P'_u$ such as in Lem\-ma~\ref{st1}.
The next goal is to prove the following result:

\ble
\lam{st2}
There exist a number\/ $r\in\dN$ and a system of\/ $\is11$
sets\/ $\pu\ne P''_u\sq P'_u$ satisfying\/
$P''_{u}\enp{\npi[u,v]}{(\pi\res n)\we r} P''_{v}$
for all\/ $u\yi v\in2^n.$ 
\ele
\bpf
Let $2^n=\ens{u_j}{j<K}$ be an arbitrary enumeration of all
dyadic sequences of length $n$; $K=2^n,$ of course.
The method of proof will be to define, for any $k\le K,$
a number $r_k\in\dN$ and 
a system of $\is11$ sets $\pu\ne Q^k_{u_j}\sq P'_{u_j}$,
$j<k,$ by induction
on $k$ so that
\ben
\fenu
\itla{**}\msur
$Q^k_{u_i}\enp{\npi[u_i,u_j]}{(\pi\res n)\we r_k} Q^k_{u_j}$ 
for all $i<j<k$. \
(Where $(\pi\res n)\we r$ is the extension of the finite sequence
$\pi\res n$ by $r$ as the new rightmost term.)
\een
After this is done, $r=r_K$ and the sets
$P''_u=Q^K_u$ prove the lemma.

We begin with $k=2.$
Then $P'_{u_0}\enq{\npi[u_0,u_1]}\pi n P'_{u_1}$ 
by \ref{z'}, and hence there exist points 
$\ang{x_0,\xi_0}\in P'_{u_0}$,   
$\ang{x_1,\xi_1}\in P'_{u_1}$ such that 
$\ang{x_0,\xi_0}\enq{\npi[u_0,u_1]}\pi n\ang{x_1,\xi_1}.$ 
Then $\xi_0\Et\xi_1,$ so that there is a number $r\in\dN$ 
with $\xi_0(n)\sd\xi_1(n)\sq r_2$.
Note that for any $p\in\dN$ and any points
$\ang{x,\xi}\yd\ang{y,\eta}\in\rnd,$
$\ang{x,\xi}\enp{\npi[u_0,u_1]}{(\pi\res n)\we r}\ang{y,\eta}$
is equivalent to
the conjunction
\dm
\ang{x,\xi}\enq{\npi[u_0,u_1]}\pi n\ang{y,\eta}
\;\;\land\;\;
\xi(n)\sd\eta(n)\sq r\,.
\dm
It follows that the sets
\dm
\bay{rcl}
S_{0}&\!\!=\!\!&
\ens{\ang{x,\xi}\in P'_{u_0}}
{\sus\ang{y,\eta}\in P'_{u_1}\,
\big(
\ang{x,\xi}\enp{\npi[u_0,u_1]}{(\pi\res n)\we r}\ang{y,\eta}
\big)}\,,\quad\text{and}\\[1\dxii]
S_{1}&\!\!=\!\!&
\ens{\ang{y,\eta}\in P'_{u_1}}
{\sus\ang{x,\xi}\in P'_{u_0}\,
\big(
\ang{x,\xi}\enp{\npi[u_0,u_1]}{(\pi\res n)\we r}\ang{y,\eta}
\big)}
\eay
\dm
are $\is11$ and non-empty 
(contain resp.\ $\ang{x_0,\xi_0}$ and $\ang{x_1,\xi_1}$), 
and they obviously satisfy 
$S_{0}\enp{\npi[u_0,u_1]}{(\pi\res n)\we r} S_{1}$.
Therefore by Corollary~\ref{suz2} there exists a 
system of $\is11$ sets $\pu\ne Q^2_u\sq P'_u\yt u\in2^n,$ 
such that $Q^2_{u_0}=S_{0}$, $Q^2_{u_1}=S_{1}$, 
\ref{z'} still holds, and in addition 
$Q^2_{u_0}\enp{\npi[u_0,u_1]}{(\pi\res n)\we r_2} Q^2_{u_1}$.
Put $r_2=r$.

Now let us carry out the step $k\mto k+1.$
Suppose that $r_k$ and sets $Q^k_{u_j}$, $j<k,$ satisfy
\ref{**}. 
Of all numbers $\npi[u_{j},u_k]\yt j<k,$ consider the least one.
Let this be, say, $\npi[u_{\ell},u_k],$ so that $\ell<k$ and 
$\npi[u_{\ell},u_k]\le\npi[u_{j},u_k]$ for all $j<k.$
As above there exists a number $r$
and a pair of non-empty $\is11$ sets
$S_{\ell}\sq Q^k_{u_\ell}$ and $S_{k}\sq Q^k_{u_k}$ such that
$S_{\ell}\enp{\npi[u_\ell,u_k]}{(\pi\res n)\we r} S_{k}$.
We can assume that $r\ge r_k$.
Put
\dm
Q'_{u_j}=\ens{\ang{y,\eta}\in S_{u_j}}
{\sus\ang{x,\xi}\in S_{\ell}\,
\big(
\ang{x,\xi}\enp{\npi[u_\ell,u_j]}{(\pi\res n)\we r}\ang{y,\eta}
\big)}
\dm
for all $j<k.$
The proof of Lemma~\ref{suz} shows that $Q'_{u_j}$ are non-empty
$\is11$ sets still satisfying \ref{**} in the form of
$Q'_{u_i}\enp{\npi[u_i,u_j]}{(\pi\res n)\we r} Q'_{u_j}$ 
for $i<j<k$ --- since $r\ge r_k$, and
obviously $Q'_{u_\ell}=S_\ell$.
In addition, put $Q'_{u_k}=S_k$.
Then still
$Q'_{u_\ell}\enp{\npi[u_\ell,u_k]}{(\pi\res n)\we r} Q'_{u_k}$
by the choice of $S_\ell$ and $S_k$.
We claim that also
\bus
\label{u8}
Q'_{u_j}\enp{\npi[u_j,u_k]}{(\pi\res n)\we r} Q'_{u_k}
\quad\text{for all }\,j<k\,.
\eus
Indeed we have
$Q'_{u_j}\enp{\npi[u_j,u_\ell]}{(\pi\res n)\we r} Q'_{u_\ell}$
and
$Q'_{u_\ell}\enp{\npi[u_\ell,u_k]}{(\pi\res n)\we r} Q'_{u_k}$ 
by the above.
It follows that
$Q'_{u_j}\enp p{(\pi\res n)\we r} Q'_{u_k}$,
where
$p=\tmax\ans{\npi[u_j,u_\ell],\npi[u_\ell,u_k]}.$
Thus it remains to show that $p\le\npi[u_j,u_k].$
That $\npi[u_\ell,u_k]\le\npi[u_j,u_k]$ holds by the choice
of $\ell.$
Prove that $\npi[u_j,u_\ell]\le\npi[u_j,u_k].$
Indeed in any case
\dm
\npi[u_j,u_\ell]\le\tmax\ans{\npi[u_j,u_k],\npi[u_\ell,u_k]}.
\dm
But once again $\npi[u_\ell,u_k]\le\npi[u_j,u_k],$ so 
$\npi[u_j,u_\ell]\le\npi[u_j,u_k]$ as required.

Thus \eqref{u8} is established.
It follows that
$Q'_{u_i}\enp{\npi[u_i,u_j]}{(\pi\res n)\we r} Q'_{u_j}$ 
for all $i<j\le k.$ 
We end the inductive step of the lemma
by putting $r_{k+1}=r$.
\epF{Lemma}

\ppu
{Step C: splitting to the next level}

We work with the number $r$ and sets $P''_u$ such as in
Lemma~\ref{st2}.
Put $\pi(n)=r.$
(Recall that $\vpi(n)$ was defined at Step A.) 
The next step is to split each one of the sets $P''_u$ in
order to define sets $P_{u\we i}\,\yt u\we i\in2^{n+1},$
of the next splitting level.

To begin with, put $Q_{u\we i}=P''_u$ for all $u\in2^n$ and
$i=0,1.$
It is easy to verify that 
the system of sets $Q_{u\we i}\,\yt u\we i\in2^{n+1},$
satisfies conditions \ref{z1} -- \ref{zlast} for the level
$n+1,$ except for \ref{z8} and \ref{z5}. 
In particular, \ref{z1'} was fixed at Step A, and \ref{z'}
in the form that
$Q_{u\we i}\enq{\npi[u\we i\,,\,v\we j]}{\pi}{(n+1)}Q_{v\we j}$
for all $u\we i$ and $v\we j$ in $2^{n+1}$ 
(and then \ref{z4} as well)
at Step B --- because ${(\pi\res n)}\we r={\pi\res{(n+1)}}.$ 

Recall that by definition all sets involved have no common
point with $\bigcup_kS^k_{\vpi(n)}$ by \ref{z1'}.
Therefore Corollary~\ref{ldc} is applicable.
We conclude that there exists a system of non-empty $\is11$
sets $W_{u\we i}\sq Q_{u\we i}\yt u\we i\in2^{n+1},$ 
still satisfying \ref{z'}, and also satisfying \ref{z5}. 

\ppu
{Step D: genericity}

We have to further shrink the sets
$W_{u\we i}\yt u\we i \in 2^{n+1},$
obtained at Step C, in order to satisfy
\ref{z8}, the last condition not yet fulfilled in the
course of the construction.
The goal is to define a new system of $\is11$ sets
$\pu\ne P_{u\we i}\sq W_{u\we i}\,\yt u\we i\in2^{n+1},$
such that still \ref{z'} holds, and in addition
$P_{u\we i}\in D_n$ for all $u\we i\in2^{n+1},$ where
$D_n$ is the \dd nth open dense subset of $\dP$ coded in $\mm$.

Take any $u_0\we i_0\in2^{n+1}.$
As $D_n$ is a dense subset of $\dP,$ there exists a
set $W_0\in D_n,$ therefore, a non-empty $\is11$ set, such that
$W_0\sq W_{u_0\we i_0}$.
It follows from Lemma~\ref{suz} that there exists a system of
non-empty $\is11$
sets $W'_{u\we i}\sq W_{u\we i}\yt u\we i\in2^{n+1},$ 
still satisfying \ref{z'}, and such that
$W'_{u_0\we i_0}=Q_0$.

Now take any other $u_1\we i_1\ne u_0\we i_0$ in $2^{n+1}.$
The same construction yields a system of non-empty $\is11$
sets $W''_{u\we i}\sq W'_{u\we i}\yt u\we i\in2^{n+1},$ 
still satisfying \ref{z'}, and such that
$W''_{u_1\we i_1}=W_1\sq W'_{u_1\we i_1}$ is a set in $D_n$.

Iterating this construction $2^{n+1}$ times, we obtain a
system of sets $P_{u\we i}$ satisfying \ref{z8} as well as
all other conditions in the list \ref{z1} -- \ref{zlast},
as required.\vtm

\qeDD{Construction and Case~2 of Theorem~\ref{mt'}}\vtm

\qeDD{Theorems~\ref{mt'} and \ref{mt}}\vtm

\bibliographystyle{plain}

\end{document}